\documentclass[12pt]{article}
\usepackage{amsmath}
\usepackage{amsthm}
\usepackage{amsfonts}
\usepackage{amssymb}
\usepackage{graphicx}
\usepackage{xypic}
\usepackage{a4wide}
\usepackage{enumerate}
\usepackage{mathrsfs}
\newtheorem{theorem}{Theorem}[section]
\newtheorem{proposition}[theorem]{Proposition}
\newtheorem{corollary}[theorem]{Corollary}
\newtheorem{lemma}[theorem]{Lemma}
\theoremstyle{definition}

\newtheorem{definition}[theorem]{Definition}
\newtheorem{example}[theorem]{Example}

\newtheorem{remark}[theorem]{Remark}
\newtheorem{remarks}[theorem]{Remarks}

\newcommand{\bcomod}[2]{{}^{#1}\mathcal{M}^{#2}}
\newcommand{\cat}[1]{\mathbf{#1}}
\newcommand{\cohom}[3]{\mathrm{h}_{#1}(#2,#3)}
\newcommand{\coring}[1]{\mathfrak{#1}}
\newcommand{\cotensor}[1]{\square_{#1}}
\renewcommand{\hom}[3]{\mathrm{Hom}_{#1}(#2,#3)}
\newcommand{\lcomod}[1]{{}^{#1}\mathcal{M}}
\newcommand{\ldual}[1]{{ }^*#1}
\newcommand{\rdual}[1]{#1^*}
\newcommand{\lDual}{{ }^*(-)}
\newcommand{\rDual}{(-)^*}
\newcommand{\rcomod}[1]{\mathcal{M}^{#1}}
\newcommand{\rmod}[1]{\mathcal{M}_{#1}}
\newcommand{\lmod}[1]{{}_{#1}\mathcal{M}}
\newcommand{\tensor}[1]{\otimes_{#1}}
\newcommand{\nat}[2]{\mathrm{Nat}(#1,#2)}
\newcommand{\h}[2]{\mathrm{H}(#1,#2)}

\date{\empty}

\begin{document}

\title{Adjoint and Frobenius Pairs of Functors for Corings}

\author{M.~Zarouali-Darkaoui\\
\normalsize Departamento de \'{A}lgebra \\
\normalsize Facultad de Ciencias \\
\normalsize Universidad de Granada\\
\normalsize E18071 Granada, Spain \\
\normalsize E-mail: \textsf{zaroual@correo.ugr.es}}

\maketitle

\section*{Introduction}

Corings were introduced by M.~Sweedler in \cite{Sweedler:1975} as a
generalization of coalgebras over commutative rings to the case of
non-commutative rings, to give a formulation of a predual of the
Jacobson-Bourbaki's theorem for intermediate extensions of division
ring extensions. Thus, a coring over an associative ring with unit
$A$ is a comonoid in the monoidal category of all $A$-bimodules.
Recently, motivated by an observation of M. Takeuchi, namely that an
entwining structure (resp. an entwined module) can be viewed as a
suitable coring (resp. as a comodule over a suitable coring),
T.~Brzezi\'{n}ski, has given in \cite{Brzezinski:2002} some new
examples and general properties of corings. Among them, a study of
Frobenius corings is developed, extending previous results on
entwining structures \cite{Brzezinski:2000} and relative Hopf
modules \cite{Caenepeel/Militaru/Zhu:1997}.

A pair of functors $(F,G)$ is said to be a \emph{Frobenius pair}
\cite{Castano/Gomez/Nastasescu:1999} if $G$ is at the same time a
left and right adjoint to $F$. That is a standard name which we use
instead of Morita's original ``\emph{strongly adjoint pairs}''
\cite{Morita:1965}. The functors $F$ and $G$ are known as
\emph{Frobenius functors} \cite{Caenepeel/Militaru/Zhu:1997}. The
study of Frobenius functors was motivated by a paper of K. Morita,
where he proved \cite[Theorem 5.1]{Morita:1965} that given a ring
extension $i:A\rightarrow B$, the induction functor
$-\tensor{A}B:\rmod{A}\rightarrow \rmod{B}$ is a Frobenius functor
if and only if the morphism $i$ is Frobenius in the sense of
\cite{Kasch:1961} (see also \cite{Nakayama/Tsuzuku:1960}). The dual
result for coalgebras over fields was proved in \cite[Theorem
3.5]{Castano/Gomez/Nastasescu:1999} and it states that the
co-restriction functor $(-)_\varphi:\mathcal{M}^{C}\rightarrow
\mathcal{M}^{D}$ associated to a morphism of coalgebras
$\varphi:C\rightarrow D$ is Frobenius if and only if $C_D$ is
quasi-finite and injective and there exists an isomorphism of
bicomodules $\cohom{D}{{}_CC_D}{D} \cong {}_DC_C$ (here,
$\cohom{D}{C}{-}$ denotes the ``cohom'' functor). Since corings
generalize both rings and coalgebras over fields, it is natural to
guess that \cite[Theorem 5.1]{Morita:1965} and \cite[Theorem
3.5]{Castano/Gomez/Nastasescu:1999} are specializations of a general
statement on homomorphisms of corings. In this paper, we find such a
result (Theorem \ref{23}), and we introduce the notion of a (right)
Frobenius extension of corings (see Definition \ref{Frobenius
extension}). To prove such a result, we study adjoint pairs and
Frobenius pairs of functors between categories of comodules over
rather general corings. Precedents for categories of modules, and
categories of comodules over coalgebras over fields are contained in
\cite{Morita:1965}, \cite{Takeuchi:1977},
\cite{Castano/Gomez/Nastasescu:1999} and
\cite{Caenepeel/DeGroot/Militaru:2002}. More recently,
\emph{Frobenius corings} (i.e., corings for which the functor
forgetting the coaction is Frobenius), have been intensively studied
in
\cite{Brzezinski:2002,Brzezinski:2003,Brzezinski/Gomez:2003,Brzezinski/Wisbauer:2003,Caenepeel/DeGroot/Militaru:2002}.

We think that our general approach produces results of independent
interest, beyond the aforementioned extension  to the coring
setting of \cite[Theorem 5.1]{Morita:1965} and \cite[Theorem
3.5]{Castano/Gomez/Nastasescu:1999}, and contributes to the
understanding of the subtle behavior of the cotensor product
functor for corings. In fact, our general results, although
sometimes rather technical, have  other applications and will
probably find more. For instance, the author has used them
to prove new results on equivalences of comodule categories over
corings \cite{Zarouali:2004}. Moreover, when applied to corings
stemming from different algebraic theories of current interest,
they boil down to new (more concrete) results. As an illustration,
we consider entwined modules over an entwining structure in
Section \ref{entwined}, and graded modules over $G$-sets in
Section \ref{graded}.

The paper is organized as follows. After Section \ref{notations},
devoted to fix some basic notations, Section
\ref{Frobeniusgeneral} deals with adjoint and Frobenius pairs on
categories of comodules. Some refinements of results from
\cite{Gomez:2002} and \cite[\S 23]{Brzezinski/Wisbauer:2003} on
the representation as cotensor product functors of certain
functors between comodule categories are needed and are thus
included in Section \ref{Frobeniusgeneral}. From our general
discussion on adjoint pairs of cotensor product functors we will
derive our main general result on Frobenius pairs between comodule
categories (Theorem \ref{20a}) that extends the known
characterizations in the setting of modules over rings and of
comodules over coalgebras. In the first case, the key property to
derive the result on modules from Theorem \ref{20a} is the
separability of the trivial corings (see Remark \ref{damodulos}).
In the case of coalgebras, the fundamental additional property is
the duality between finite left and right comodules. We already
consider a much more general situation in Section
\ref{condualidad}, where we introduce the class of so called
\emph{corings having a duality} for which we prove
characterizations of Frobenius pairs that are similar to the
coalgebra case.

Section \ref{Frobmor} is the leitmotiv of the paper. After the
technical development of sections \ref{Frobeniusgeneral} and
\ref{condualidad}, our main results follow without difficulty. We
prove in particular that the induction functor 
$-\otimes_{A}B:\mathcal{M}^{\coring{C}}
\rightarrow\mathcal{M}^{\coring{D}}$ associated to a homomorphism
$(\varphi,\rho):\coring{C}\rightarrow\coring{D}$ of corings 
$\coring{C}$ and
$\coring{D}$ flat over their respective base rings $A$ and $B$ is
Frobenius if and only if the $\coring{C} -\coring{D}$-bicomodule
$\coring{C}\otimes_{A}B$ is quasi-finite and injector as a right
$\coring{D}$-comodule and there exists an isomorphism of
$\coring{D}-\coring{C}$-bicomodules
$\cohom{\coring{D}}{\coring{C}\otimes_{A}B}{\coring{D}}\simeq
B\otimes _{A}\coring{C}$ (Theorem \ref{23}). We show as well how
this theorem unifies previous results for ring homomorphisms
\cite[Theorem 5.1]{Morita:1965}, coalgebra maps \cite[Theorem
3.5]{Castano/Gomez/Nastasescu:1999}, and Frobenius corings
\cite[27.10, 28.8]{Brzezinski/Wisbauer:2003}.

In Section \ref{entwined}, we specialize one of the general
results on corings to entwining structures.

In Section \ref{graded}, we particularize our results in the
previous sections to the coring associated to a $G$-graded algebra
and a $G$-set, where $G$ ia a group. We obtain a series of new
results for graded modules by $G$-sets.

\section{Basic notations}\label{notations}

Throughout this paper and unless otherwise stated, $k$ denote a
commutative ring (with unit), $A,$ $A',$ $A'',$ and $B$ denote
associative and unitary algebras over $k$, and $\coring{C},$
$\coring{C}',$ $\coring{C}'',$ and $\coring{D}$ denote corings over
$A,$ $A',$ $A'',$ and $B$, respectively. We recall from
\cite{Sweedler:1975} that an $A$-\emph{coring} consists of an
$A$-bimodule $\coring{C}$ with two $A$-bimodule maps
\[
\Delta : \coring{C} \longrightarrow \coring{C} \tensor{A}
\coring{C}, \qquad \epsilon : \coring{C} \longrightarrow A
\]
such that $(\coring{C} \tensor{A} \Delta) \circ \Delta = (\Delta
\tensor{A} \coring{C}) \circ \Delta$ and $(\epsilon \tensor{A}
\coring{C}) \circ \Delta = ( \coring{C} \tensor{A} \epsilon) \circ
\Delta = 1_\coring{C}$. A \emph{right} $\coring{C}$-\emph{comodule}
is a pair $(M,\rho_M)$ consisting of a right $A$-module $M$ and an
$A$-linear map $\rho_M: M \rightarrow M \tensor{A} \coring{C}$ (the
coaction) satisfying $(M \tensor{A} \Delta) \circ \rho_M = (\rho_M
\tensor{A} \coring{C}) \circ \rho_M$, $(M \tensor{A} \epsilon) \circ
\rho_M = 1_M$. A \emph{morphism} of right $\coring{C}$-comodules
$(M,\rho_M)$ and $(N,\rho_N)$ is a right $A$-linear map $f: M
\rightarrow N$ such that $(f \tensor{A} \coring{C}) \circ \rho_M =
\rho_N \circ f$; the $k$-module of all such morphisms will be
denoted by $\hom{\coring{C}}{M}{N}$. The right
$\coring{C}$-comodules together with their morphisms form the
additive category $\rcomod{\coring{C}}$. Coproducts and cokernels
(and then inductive limits) in $\rcomod{\coring{C}}$ exist and they
coincide respectively with coproducts and cokernels in the category
of right $A$-modules $\rmod{A}$. If ${}_A\coring{C}$ is flat, then
$\rcomod{\coring{C}}$ is a Grothendieck category.
The converse is not true in general (see \cite[Example 1.1]
{ElKaoutit/Gomez/Lobillo:2001unp}). When $\coring{C} = A$
with the trivial $A$-coring structure, then $\rcomod{A} =
\rmod{A}$.

Now assume that the $A'-A$-bimodule $M$ is also a left comodule
over an $A'$-coring $\coring{C}'$ with structure map $\lambda_M :
M \rightarrow \coring{C}' \tensor{A'} M$. Assume moreover that $\rho_M$
is $A'$-linear, and $\lambda_M$ is $A$-linear. It is clear that $\rho_M
: M \rightarrow M \tensor{A} \coring{C}$ is a morphism of left
$\coring{C}'$-comodules if and only if $\lambda_M : M \rightarrow
\coring{C}' \tensor{A'} M$ is a morphism of right
$\coring{C}$-comodules. In this case, we say that $M$ is a
$\coring{C}'-\coring{C}$-bicomodule. The
$\coring{C}'-\coring{C}$-bicomodules are the objects of a category
$\bcomod{\coring{C}'}{\coring{C}}$, whose morphisms are defined in
the obvious way.

Let $Z$ be a left $A$-module and $f: X \rightarrow Y$ a morphism
in $\rmod{A}$. Following \cite[40.13]{Brzezinski/Wisbauer:2003} we
say that $f$ is \emph{$Z$-pure} when the functor $-\tensor{A} Z$
preserves the kernel of $f$. If $f$ is $Z$-pure for every $Z \in
\lmod{A}$ then we say simply that $f$ is \emph{pure} in
$\rmod{A}$.  The notation $\otimes$ will stand for the tensor
product over $k$.

\section{Frobenius functors between categories of
comodules}\label{Frobeniusgeneral}

Let $T$ be a $k$-algebra, and $M\in{}^{T}\mathcal{M}^{A}$. Let
$\varphi:T\rightarrow \operatorname{End}_A(M)$ the morphism of
$k$-algebras given by the right $T$-module structure of the bimodule
$_TM_A$. Now, suppose moreover that $M\in\mathcal{M}^{\coring{C}}$.
Then $\operatorname{End}_{\coring{C}}(M)$ is a subalgebra of
$\operatorname{End}_A(M).$ We have that $\varphi(T)\subset
\operatorname{End}_{\coring{C}}(M)$ if and only if $\rho_M$ is
$T$-linear.
Hence the left $T$-module structure of a $T-\coring{C}$-bicomodule
$M$ can be described as a morphism of $k$-algebras $
\varphi:T\rightarrow \operatorname{End}_{\coring{C}}(M)$. Given a
$k$-linear functor $F:\mathcal{M}^{\coring{C}}\rightarrow\mathcal{M}
^{\coring{D}}$, and $M\in{}^{T}\mathcal{M}%
^{\coring{C}}$, the algebra morphism $\xymatrix{T \ar[r]^-{\varphi}
& \operatorname{End}_{\coring{C}}(M) \ar[r]^-{F(-)} &\operatorname
{End}_{\coring{D}}(F(M))}$ defines a $T-\coring{D}$-bicomodule
structure on $F(M)$. We have then two $k$-linear bifunctors
\[
-\otimes_TF(-),\,F(-\otimes_T-):\mathcal{M}%
^{T}\times{}^{T}\mathcal{M}^{\coring{C}}\rightarrow\mathcal{M}%
^{\coring{D}}.
\]
Let $\Upsilon_{T,M}$ be the unique isomorphism of
$\coring{D}$-comodules making the following diagram commutative

\begin{equation}\label{T}
\xymatrix{  T\otimes_TF(M)\ar[rd]_\simeq
\ar[rr]^{\Upsilon_{T,M}} & &  F(T\otimes_TM)\ar[ld]^\simeq \\
 & F(M) & }
\end{equation}
for every $M\in{}^{T}\mathcal{M}^{\coring{C}}.$ We have that
$\Upsilon_{T,M}$ is natural in $T.$ By the theorem of Mitchell
\cite[Theorem 3.6.5]{Popescu:1973}, there exists a unique natural
transformation
\[
\Upsilon_{-,M}:-\otimes_TF(M)\rightarrow F(-\otimes_TM)
\]
extending the natural transformation $\Upsilon_{T,M}.$ We refer to
\cite{Gomez:2002} for more details.

\begin{remark}\label{Mitchell}
Mitchell's \cite[Theorem 3.6.5]{Popescu:1973} holds also if we
only suppose that the target category $\mathscr{C'}$ is
preadditive and has coproducts, or if the category $\mathscr{C'}$
is preadditive and the functor $S$ preserves coproducts (notations
as in \cite[Theorem 3.6.5]{Popescu:1973}). This fact is used to
show that the natural transformation $\Upsilon$ exists for every
$k$-linear functor
$F:\mathcal{M}^{\coring{C}}\rightarrow\mathcal{M}^{\coring{D}}$
even if the category $\mathcal{M}^{\coring{D}}$ is not abelian.
Note also that its corollary \cite[Corollary 3.6.6]{Popescu:1973}
also holds if we suppose only that the category $\mathscr{C'}$ is
preadditive.
\end{remark}

Let $M\in{}\bcomod{\coring{C'}}{\coring{C}}$ and
$N\in{}^{\coring{C}}\mathcal{M}^{\coring{C}''}$. The map
\[
\omega_{M,N}=\rho_M\otimes_AN-M\otimes_A\lambda_N:M\otimes_AN
\rightarrow M\otimes_A\coring{C}\otimes_AN
\]
is a $\coring{C}'-\coring{C}^{\prime\prime}$-bicomodule map. Its kernel in $_{A'}\mathcal{M}_{A''}$ is the \emph{cotensor product} of $M$ and $N$, and it is denoted by $M\square_{\coring{C}}N$.
If $\omega_{M,N}$ is $\coring{C}'_{A'}$-pure and
$_{A''}\coring{C}''$-pure, and the following
\begin{equation}
\operatorname{ker} (\omega_{M,N})\otimes_{A''}\coring{C}''
\otimes_{A''}\coring{C}'', \quad \coring{C}' \otimes_{A'}\coring{C}'
\otimes_{A'}\operatorname{ker} (\omega_{M,N}) \quad \textrm{and}
\quad \coring{C}'\otimes_{A'}\operatorname{ker}
(\omega_{M,N})\otimes_{A''}\coring{C}''
\end{equation}
are injective maps, then $M\square_{\coring{C}}N$ is the kernel of
$\omega_{M,N}$ in ${}^{\coring{C}'}\mathcal{M}^{\coring{C}''}$. This
is the case if $\omega_{M,N}$ is
$(\coring{C}'\otimes_{A'}\coring{C}')_{A'}$-pure,
$_{A''}(\coring{C}'' \otimes_{A''}\coring{C}'')$-pure, and
$\coring{C}'\otimes_{A'}\omega_{M,N}$ is $_{A''}\coring{C}''$-pure
(e.g. if $\coring{C}'_{A'}$ and $_{A''}\coring{C}''$ are flat, or if
$\coring{C}$ is a coseparable $A$-coring).

If for every $M\in{}\bcomod{\coring{C'}}{\coring{C}}$ and
$N\in{}^{\coring{C}}\mathcal{M}^{\coring{C}''}$, $\omega_{M,N}$ is
$\coring{C}'_{A'}$-pure and $_{A''}\coring{C}''$-pure, then we have
a $k$-linear bifunctor
\begin{equation}\label{cotensorbifunctor}
\xymatrix{-\cotensor{\coring{C}}-:{}^{\coring{C}'}\mathcal{M}^{\coring{C}}
\times{}^{\coring{C}}\mathcal{M}^{\coring{C}''}\ar[r]&
{}^{\coring{C}'}\mathcal{M}^{\coring{C}''}}.
\end{equation} If in particular $\coring{C}'_{A'}$ and
$_{A''}\coring{C}'' $ are flat, or if $\coring{C}$ is a coseparable
$A$-coring, then the bifunctor \eqref{cotensorbifunctor} is well
defined.

By a proof similar to that of \cite[II.1.3]{AlTakhman:1999}, we
have, for every $M\in\mathcal{M}^{\coring{C}}$, that the functor
$M\square_{\coring{C}}-$ preserves direct limits.

The following lemma was used implicitly in the proof of
\cite[Proposition 3.4]{Gomez:2002}, and it will be useful for us
in the proof of the next theorem.

\begin{lemma}\label{2}
If$M\in{}\bcomod{\coring{C'}}{\coring{C}}$, and
$F:\mathcal{M}^{\coring{C}}\rightarrow\mathcal{M}^{\coring{D}}$ is a
$M$-compatible $k$-linear functor in the sense of \cite[p.
210]{Gomez:2002}, which preserves
coproducts, then for all $X\in\mathcal{M}^{A'},$%
\[
\Upsilon_{X\otimes_{A'}\coring{C}',M}(X\otimes_{A'%
}\lambda_{F(M)})=F(X\otimes_{A'}\lambda_M)\Upsilon_{X,M}.
\]
\end{lemma}

\begin{proof}
Let us consider the diagram
\[
\xymatrix{
 X\otimes_{A'}F(M) \ar[rrrr]^{X\otimes_{A'}\lambda
_{F(M)}} \ar[rrd]^{X\otimes_{A'}F(\lambda_M)}
\ar[dd]_{\Upsilon_{X,M}}& & & & X \otimes_{A'}\coring{C}^{\prime
}\otimes_{A'}F(M) \ar[lld]^{X\otimes_{A'}\Upsilon_{\coring{C}',M}}
\ar[dd]^{\Upsilon_{X\otimes
_{A'}\coring{C}',M}}\\
& & X\otimes_{A'}F(\coring{C} '\otimes_{A'}M)
\ar[rrd]^{\Upsilon_{X,\coring{C}'\otimes_{A'}M}} & & \\
F(X\otimes_{A'}M) \ar[rrrr]_{F(X\otimes_{A'}\lambda _M)}& &  & &
F(X\otimes_{A'}\coring{C}^{\prime }\otimes_{A'}M) }
\]
The commutativity of the top triangle follows from the definition
of $\lambda_{F(M)}$, while the right triangle commutes by
\cite[Lemma 3.3]{Gomez:2002} (we take
$S=T=A'$, and $Y=\coring{C}'%
$), and the left triangle is commutative since $\Upsilon_{X,-}$
is natural. Therefore, the commutativity of the rectangle holds.
\end{proof}

 A closer analysis of \cite[Theorem 3.5]{Gomez:2002} gives the
following generalization of \cite[Proposition 2.1]{Takeuchi:1977}
and \cite[23.1(1)]{Brzezinski/Wisbauer:2003}. Recall from
\cite{Guzman:1989} that a coring $\coring{C}$ is said to be
\emph{coseparable} if the comultiplication map
$\Delta_{\coring{C}}$ is a split monomorphism of
$\coring{C}$-bicomodules. Of course, the trivial $A$-coring
$\coring{C} = A$ is coseparable and, henceforth, every result for
comodules over coseparable corings applies in particular for
modules over rings.

\begin{theorem}\label{3}
Let
$F:\mathcal{M}^{\coring{C}}\rightarrow\mathcal{M}^{\coring{D}}$ be
a $k$-linear functor, such that
\begin{enumerate}[(I)]
\item $_B\coring{D}$ is flat and $F$ preserves the kernel of
$\rho_N\otimes_A\coring{C}-N\otimes_A%
\Delta_{\coring{C}}$ for every $N \in \rcomod{\coring{C}}$, or
\item $\coring{C}$ is a coseparable $A$-coring and the categories
$\mathcal{M}^{\coring{C}}$ and $\mathcal{M}^{\coring{D}}$ are
abelian.
\end{enumerate}
Assume that at least one of the following statements holds
\begin{enumerate}
\item $\coring{C}_A$ is projective, $F$ preserves coproducts,
and $\Upsilon_{N,\coring{C}}$, $\Upsilon_{N\otimes_A\coring{C}%
,\coring{C}}$ are isomorphisms for all
$N\in\mathcal{M}^{\coring{C}}$ (e.g. , if $A$ is semisimple and $F$
preserves coproducts), or \item  $\coring{C}_A$ is flat, $F$
preserves direct limits, and $\Upsilon _{N,\coring{C}}$,
$\Upsilon_{N\otimes_A\coring{C},\coring{C}}$ are isomorphisms for
all $N\in\mathcal{M}^{\coring{C}}$ (e.g. , if $A$ is a von Neumann
regular ring and $F$ preserves direct limits), or \item $F$
preserves inductive limits (e.g., if $F$ has a right adjoint).
\end{enumerate}
 Then $F$ is naturally equivalent to
$-\square_{\coring{C}}F(\coring{C}).$
\end{theorem}

\begin{proof}
At first, note that if $\coring{C}_A$ is projective, then the right
$A$-module $\coring{C}\otimes_A\coring{C}$ is projective (by
\cite[Example 3 p. 105, Proposition VI.9.5]{Stenstrom:1975}). Hence,
if $\coring{C}_A$ is projective and $F$ preserves coproducts, then
$F$ is $M$-compatible in the sense of \cite[p. 210]{Gomez:2002}, for
all
$M\in{}^{\coring{C}}\mathcal{M}^{\coring{C}%
}.$ In each case, we have $F$ is $\coring{C}$-compatible where
$\coring{C}\in{}^{\coring{C}}\mathcal{M}^{\coring{C}}.$ Therefore,
by \cite[Proposition 3.4]{Gomez:2002}, $F(\coring{C})$ can be
viewed as a $\coring{C}-\coring{D}$-bicomodule. From Lemma
\ref{2}, and since $\Upsilon_{-,\coring{C}}$ is a natural
transformation, we have, for every $N\in
\mathcal{M}^{\coring{C}}$, the commutativity of the following
diagram with
exact rows in $\mathcal{M}^{\coring{D}}$%
\[%
\xymatrix{
 0 \ar[r] & N\square_{\coring{C}}F(\coring{C})\ar[r] &
 N\otimes_AF(\coring{C})
\ar[d]^\simeq_{\Upsilon_{N,\coring{C}}}
\ar[rrr]^-{\rho_N\otimes_AF(\coring{C})
-N\otimes_A%
\lambda_{F(\coring{C})}} &  & & N\otimes
_A\coring{C}\otimes_AF(\coring{C})\ar[d]^{\simeq
}_{\Upsilon_{N\otimes_A\coring{C},\coring{C}}} \\
 0 \ar[r] &  F(N)\ar[r]^-{F(\rho
_N)} & F(N\otimes_A\coring{C})
\ar[rrr]^-{F(\rho_N\otimes_A\coring{C}-N\otimes_A%
\Delta_{\coring{C}})} & & & F(N\otimes
_A\coring{C}\otimes_A\coring{C}).}
\]
The exactness of the bottom sequence is assumed in the case $(I)$.
For the case $(II)$, it follows by factorizing the map
$\omega_{N,\coring{C}} =
\rho_N\otimes_A\coring{C}-N\otimes_A%
\Delta_{\coring{C}}$ through its image, and using the facts that the
sequence $\xymatrix{0 \ar[r] &  N\ar[r]^-{\rho _N} & N
\otimes_A\coring{C} \ar[r]^-{\omega_{N,\coring{C}}}  & N \otimes
_A\coring{C}\otimes_A\coring{C}}$ is split exact in
$\rcomod{\coring{C}}$ in the sense of
\cite[40.5]{Brzezinski/Wisbauer:2003}, and that additive functors
between abelian categories preserve split exactness. By the
universal property of kernels, there exists a unique isomorphism
$\eta_N:N\square _{\coring{C}}F(\coring{C}) \rightarrow F(N)$ in
$\mathcal{M}^{\coring{D}}$ making commutative the above diagram. It
easy to show that $\eta$ is natural. Hence
$F\simeq-\square_{\coring{C}}F(\coring{C}).$
\end{proof}

As an immediate consequence of the last theorem we have the
following generalization of Eilenberg-Watts Theorem
\cite[Proposition VI.10.1]{Stenstrom:1975}.

\begin{corollary}\label{4}
Let $F:\mathcal{M}^{\coring{C}}\rightarrow\mathcal{M}^{\coring{D}}$
be a $k$-linear functor.\begin{enumerate}[1.]\item If $_B\coring{D}$
is flat and $A$ is a semisimple ring (resp. a von Neumann regular
ring), then the following statements are
equivalent\begin{enumerate}[(a)]\item $F$ is left exact and
preserves coproducts (resp. left exact and preserves direct limits);
\item $F\simeq-\square_{\coring{C}}M$ for some bicomodule
$M\in{}^{\coring{C}%
}\mathcal{M}^{\coring{D}}$.\end{enumerate}
 \item If
$_A\coring{C}$ and $_B\coring{D}$ are flat, then the following
statements are equivalent \begin{enumerate}[(a)]\item $F$ is exact
and preserves inductive limits;\item $F\simeq-\square_{\coring{C}}M$
for some
bicomodule $M\in{}^{\coring{C}%
}\mathcal{M}^{\coring{D}}$ which is coflat in
$^{\coring{C}}\mathcal{M}%
$. \end{enumerate} \item If $\coring{C}$ is a coseparable
$A$-coring and the categories $\mathcal{M}^{\coring{C}}$ and
$\mathcal{M}^{\coring{D}}$ are abelian, then the following
statements are equivalent \begin{enumerate}[(a)]\item $F$
preserves inductive
limits;\item $F$ preserves cokernels and $F\simeq-\square_{\coring{C}%
}M$ for some bicomodule
$M\in{}^{\coring{C}}\mathcal{M}^{\coring{D}}$.
\end{enumerate}\item If $\coring{C}=A$ and the category
$\mathcal{M}^{\coring{D}}$ is abelian, then the following
statements are equivalent
\begin{enumerate}[(a)]\item $F$ has a right adjoint;\item $F$
preserves inductive limits;\item
$F\simeq-\otimes_AM$ for some bicomodule $M\in{}^{A}\mathcal{M}%
^{\coring{D}}.$\end{enumerate}
\end{enumerate}
\end{corollary}

A bicomodule $N \in \bcomod{\coring{C}}{\coring{D}}$ is said to be
\emph{quasi-finite} as a right $\coring{D}$-comodule if the
functor $- \tensor{A} N : \rmod{A} \rightarrow
\rcomod{\coring{D}}$ has a left adjoint $\cohom{\coring{D}}{N}{-} :
\rcomod{\coring{D}} \rightarrow \rmod{A}$, the \emph{cohom
functor}. If $\omega_{Y,N}$ is $\coring{D} \tensor{B}
\coring{D}$-pure for every right $\coring{C}$-comodule $Y$ (e.g.,
${}_B\coring{D}$ is flat or $\coring{C}$ is coseparable) then
$N_{\coring{D}}$ is quasi-finite if and only if $-
\cotensor{\coring{C}} N : \rcomod{\coring{C}} \rightarrow
\rcomod{\coring{D}}$ has a left adjoint, which we still to denote
by $\cohom{\coring{D}}{N}{-}$ \cite[Proposition 4.2]{Gomez:2002}. The
particular case of the following statement when the cohom is exact
generalizes \cite[Corollary 3.12]{AlTakhman:2002}.

\begin{corollary}\label{5}
Let $N\in{}^{\coring{C}}\mathcal{M}^{\coring{D}}$ be a bicomodule,
quasi-finite as a right $\coring{D}$-comodule, such that
$_A\coring{C}$ and $_B\coring{D}$ are flat. If the cohom functor
$\cohom{\coring{D}}{N}{-}$ is exact or if $\coring{D}$ is a
coseparable $B$-coring, then we have
\[
\cohom{\coring{D}}{N}{-}
\simeq-\square_{\coring{D}}\cohom{\coring{D}}{N}{\coring{D}}:
\mathcal{M}^{\coring{D}}\rightarrow \mathcal{M}^{\coring{C}}.
\]
\end{corollary}

\begin{proof}
The functor $\cohom{\coring{D}}{N}{-}$ is $k$-linear and preserves
inductive limits, since it is a left adjoint to the $k$-linear
functor
$-\square_{\coring{C}}N:\mathcal{M}^{\coring{C}}\rightarrow
\mathcal{M}^{\coring{D}}$ (by \cite[Proposition 4.2]{Gomez:2002}).
Hence Theorem \ref{3} achieves the proof.
\end{proof}

Now we will use the following generalization of \cite[Lemma
2.2]{Takeuchi:1977}.

\begin{lemma}\label{18}
Let $\Lambda$, $\Lambda'$ be
bicomodules in $^{\coring{D}%
}\mathcal{M}^{\coring{C}}$ and $G=-\square_{\coring{D}}\Lambda$,
$G'=-\square_{\coring{D}}\Lambda'.$ Suppose moreover that
${}_A\coring{C}$ is flat and $B$ is a von Neumann regular ring, or
${}_A\coring{C}$ is flat and $G$ and $G'$ are cokernel preserving,
or $\coring{D}$ is a coseparable coring. Then
\[
\nat{G}{G'}\simeq \hom{(\coring{D},\coring{C})}{\Lambda}{\Lambda'}.
\]
\end{lemma}

\begin{proof}
Let $\alpha:G\rightarrow G'$ be a natural transformation. By
\cite[Lemma 3.2(1)]{Gomez:2002}, $\alpha_{\coring{D}}$ is left
$B$-linear. For the rest of the proof it suffices to replace
$\otimes$ by $\otimes_B$ in the proof of \cite[Lemma
4.1]{Caenepeel/DeGroot/Militaru:2002}.
\end{proof}

The following proposition generalizes \cite[Theorem
2.1]{Morita:1965} from bimodules over rings to bicomodules over
corings.

\begin{proposition}\label{19a}
Suppose that $_A\coring{C}$, $\coring{C}_A$, $_B\coring{D}$ and
$\coring{D}_B$ are flat. Let $X\in{}^{\coring{C}}\mathcal{M}%
^{\coring{D}}$ and
$\Lambda\in{}^{\coring{D}}\mathcal{M}^{\coring{C}}$. Consider the
following properties:
\begin{enumerate}[(1)]
\item \label{adj1} $-\square _{\coring{C}}X$ is left adjoint to
$-\square_{\coring{D}}\Lambda$; \item \label{adj2} $\Lambda$ is
quasi-finite as a right $\coring{C}$-comodule and $-
\square_{\coring{C}} X \simeq \cohom{\coring{C}}{\Lambda}{-}$;
\item \label{adjsep} $\Lambda$ is quasi-finite as a right $\coring{C}%
$-comodule and $X\simeq \cohom{\coring{C}}{\Lambda}{\coring{C}}$
in $^{\coring{C}}\mathcal{M}^{\coring{D}}$; \item
\label{semicontext} there exist bicolinear maps
\[
\psi:\coring{C}\rightarrow X\square_{\coring{D}}\Lambda\text{ and
}\omega:\Lambda\square_{\coring{C}}X\rightarrow\coring{D}
\]
in $^{\coring{C}}\mathcal{M}^{\coring{C}}$ and $^{\coring{D}}
\mathcal{M}^{\coring{D}}$ respectively, such that
\begin{equation}\label{unitcounit}
(\omega\square_{\coring{D}}\Lambda)\circ(\Lambda
\square_{\coring{C}}\psi)=\Lambda\text{ and }(X\square
_{\coring{D}}\omega)\circ(\psi\square_{\coring{C}}X)=X;
\end{equation}
\item \label{adj1symm} $\Lambda \cotensor{\coring{C}}-$ is left
adjoint to $X \cotensor{\coring{D}}-$.
\end{enumerate}
Then (\ref{adj1}) and (\ref{adj2}) are equivalent, and they imply
(\ref{adjsep}). The converse is true if $\coring{C}$ is a
coseparable $A$-coring. If ${}_A X$ and ${}_B \Lambda$ are flat, and
$\omega{}_{X,\Lambda}=\rho{}_X\otimes_B\Lambda-X\otimes_A\rho{}_\Lambda$
is pure as an $A$-linear map and
$\omega{}_{\Lambda,X}=\rho{}_\Lambda\otimes_AX-\Lambda\otimes_B\rho{}_X$
is pure as a $B$-linear map (e.g. if ${}_{\coring{C}}X$ and
${}_{\coring{D}}\Lambda$ are coflat
\cite[21.5]{Brzezinski/Wisbauer:2003} or $A$ and $B$ are von Neumann
regular rings), or if $\coring{C}$ and $\coring{D}$ are coseparable,
then (\ref{semicontext}) implies (\ref{adj1}). The converse is true
if ${}_{\coring{C}}X$ and ${}_{\coring{D}}\Lambda$ are coflat, or if
$A$ and $B$ are von Neumann regular rings, or if $\coring{C}$ and
$\coring{D}$ are coseparable. Finally, if $\coring{C}$ and
$\coring{D}$ are coseparable, or if $X$ and $\Lambda$ are coflat on
both sides, or if $A, B$ are von Neumann regular rings, then
(\ref{adj1}), (\ref{semicontext}) and (\ref{adj1symm}) are
equivalent.
\end{proposition}

\begin{proof}
The equivalence between (\ref{adj1}) and (\ref{adj2}) follows from
\cite[Proposition 4.2]{Gomez:2002}. That (\ref{adj2}) implies
(\ref{adjsep}) is a consequence of \cite[Proposition
3.4]{Gomez:2002}. If $\coring{C}$ is coseparable and we assume
(\ref{adjsep}) then, by Corollary \ref{5},
$\cohom{\coring{C}}{\Lambda}{-} \simeq - \cotensor{\coring{C}}
\cohom{\coring{C}}{\Lambda}{\coring{C}} \simeq -
\cotensor{\coring{C}} X$. That (\ref{adj1}) implies
(\ref{semicontext}) follows from Lemma \ref{18} by evaluating the
unit and the counit of the adjunction at $\coring{C}$ and
$\coring{D}$, respectively. Conversely, if we put
$F=-\cotensor{\coring{C}}X$ and $G=-\cotensor{\coring{D}}\Lambda$,
we have $GF \simeq
-\square_{\coring{C}%
}(X\square_{\coring{D}}\Lambda)$ and $FG \simeq -\square
_{\coring{D}}(\Lambda\square_{\coring{C}}X)$ by \cite[Proposition
22.6]{Brzezinski/Wisbauer:2003}$. $ Define natural transformations
$$\xymatrix@1{\eta :1_{\mathcal{M}^{\coring{C}}}\ar[r]^\simeq &
-\square
_{\coring{C}}\coring{C}\ar[r]^-{-\square_{\coring{C}}\psi} & GF}$$
 and
$$ \xymatrix@1{\varepsilon:FG\ar[r]^{-\square_{\coring{D}}\omega} &
-\square_{\coring{D}}\coring{D}\ar[r]^\simeq &
1_{\mathcal{M}^{\coring{D}}},} $$ which become the unit and the
counit of an adjunction by \eqref{unitcounit}. This gives the
equivalence between (\ref{adj1}) and (\ref{semicontext}). The
equivalence between (\ref{semicontext}) and (\ref{adj1symm}) follows
by symmetry.
\end{proof}

\begin{definition}\label{7}
Following \cite{AlTakhman:2002} and \cite{Brzezinski/Wisbauer:2003},
a bicomodule $N\in{}^{\coring{C}}\mathcal{M}^{\coring{D}}$ is called
an \emph{injector} as a right $\coring{D}$-comodule if the functor
$-\otimes_AN:\mathcal{M}^{A}\rightarrow\mathcal{M}^{\coring{D}}$
preserves injective objects.
\end{definition}

\begin{proposition}\label{19}
Suppose that $_A\coring{C}$ and $_B\coring{D}$ are flat. Let
$X\in{}^{\coring{C}}\mathcal{M}%
^{\coring{D}}$ and
$\Lambda\in{}^{\coring{D}}\mathcal{M}^{\coring{C}}.$ The following
statements are equivalent
\begin{enumerate}[(i)]
\item $-\square _{\coring{C}}X$ is left adjoint to
$-\square_{\coring{D}}\Lambda$, and $- \cotensor{\coring{C}} X$ is
left exact (or ${}_AX$ is flat or ${}_{\coring{C}}X$ is coflat);
\item $\Lambda$ is quasi-finite as a right $\coring{C}$-comodule,
$- \square_{\coring{C}} X \simeq \cohom{\coring{C}}{\Lambda}{-}$,
and $- \square_{\coring{C}}X$ is left exact (or ${}_AX$ is flat or
${}_{\coring{C}}X$ is coflat); \item
 $\Lambda$ is quasi-finite and injector as a right
$\coring{C}$-comodule and $X\simeq
\cohom{\coring{C}}{\Lambda}{\coring{C}}$ in
$^{\coring{C}}\mathcal{M}^{\coring{D}}$.
\end{enumerate}
\end{proposition}

\begin{proof}
First, observe that if ${}_{\coring{C}}X$ is coflat, then ${}_AX$
is flat \cite[21.6]{Brzezinski/Wisbauer:2003}, and that if ${}_AX$
is flat, then the functor $- \cotensor{\coring{C}} X$ is left
exact. Thus, in view of Proposition \ref{19a}, it suffices if we
prove that the version of $(ii)$ with $- \cotensor{\coring{C}} X$
left exact implies $(iii)$, and this last implies the version of
$(ii)$ with ${}_{\coring{C}}X$ coflat. Assume that $-
\cotensor{\coring{C}}X \simeq \cohom{\coring{C}}{\Lambda}{-}$ with
$- \cotensor{\coring{C}}X$ left exact. By \cite[Proposition
3.4]{Gomez:2002}, $X \simeq
\cohom{\coring{C}}{\Lambda}{\coring{C}}$ in
$\bcomod{\coring{C}}{\coring{D}}$. Being a left adjoint,
$\cohom{\coring{C}}{\Lambda}{-}$ is right exact and, henceforth,
exact. By \cite[Theorem 3.2.8]{Popescu:1973},
$\Lambda_{\coring{C}}$ is an injector and we have proved $(iii)$.
Conversely, if $\Lambda_{\coring{C}}$ is a quasi-finite injector
and $X \simeq \cohom{\coring{C}}{\Lambda}{\coring{C}}$ as
bicomodules, then $ - \cotensor{\coring{C}} X \simeq -
\cotensor{\coring{C}} \cohom{\coring{C}}{\Lambda}{\coring{C}}$
and, by \cite[Theorem 3.2.8]{Popescu:1973}, we get that
$\cohom{\coring{C}}{\Lambda}{-}$ is an exact functor. By Corollary
\ref{5}, $\cohom{\coring{C}}{\Lambda}{-} \simeq -
\cotensor{\coring{C}} \cohom{\coring{C}}{\Lambda}{\coring{C}}
\simeq - \cotensor{\coring{C}} X$, and ${}_{\coring{C}}X$ is
coflat.
\end{proof}

From the foregoing propositions, it is easy to deduce our
characterization of Frobenius functors between categories of
comodules over corings.

\begin{theorem}\label{20a}
Suppose that $_A\coring{C}$ and $_B\coring{D}$ are flat. Let
$X\in{}^{\coring{C}}\mathcal{M}%
^{\coring{D}}$ and
$\Lambda\in{}^{\coring{D}}\mathcal{M}^{\coring{C}}.$ The following
statements are equivalent
\begin{enumerate}[(i)]
\item $(-\cotensor{\coring{C}} X, - \cotensor{\coring{D}}
\Lambda)$ is a Frobenius pair; \item $ - \cotensor{\coring{C}} X$
is a Frobenius functor, and $\cohom{\coring{D}}{X}{\coring{D}}
\simeq \Lambda$ as bicomodules; \item there is a Frobenius pair
$(F,G)$ for $\rcomod{\coring{C}}$ and $\rcomod{\coring{D}}$ such
that $F(\coring{C}) \simeq \Lambda$ and $G(\coring{D}) \simeq X$
as bicomodules; \item $\Lambda_{\coring{C}}, X_{\coring{D}}$ are
quasi-finite injectors,  and $X \simeq
\cohom{\coring{C}}{\Lambda}{\coring{C}}$ and $\Lambda \simeq
\cohom{\coring{D}}{X}{\coring{D}}$ as bicomodules; \item
$\Lambda_{\coring{C}}, X_{\coring{D}}$ are quasi-finite, and $-
\cotensor{\coring{C}} X \simeq \cohom{\coring{C}}{\Lambda}{-}$ and
$ - \cotensor{\coring{D}} \Lambda \simeq
\cohom{\coring{D}}{X}{-}$.
\end{enumerate}
\end{theorem}

\begin{proof}
$(i) \Leftrightarrow (ii) \Leftrightarrow (iii)$ This is obvious,
after Theorem \ref{3} and \cite[Proposition 3.4]{Gomez:2002}.

$(i) \Leftrightarrow (iv)$ Follows from Proposition \ref{19}.

$(iv) \Leftrightarrow (v)$ If $X_{\coring{D}}$ and
$\Lambda_{\coring{C}}$ are quasi-finite, then ${}_AX$ and
${}_B\Lambda$ are flat. Now, apply Proposition \ref{19}.
\end{proof}

From Proposition \ref{19a} and Proposition \ref{19} (or Theorem
\ref{20a}) we get the following

\begin{theorem}\label{leftrightFrobenius}
Let $X\in{}^{\coring{C}}\mathcal{M}%
^{\coring{D}}$ and
$\Lambda\in{}^{\coring{D}}\mathcal{M}^{\coring{C}}.$ Suppose that
$_A\coring{C}$, $\coring{C}_A$, $_B\coring{D}$ and $\coring{D}_B$
are flat. The following statements are equivalent
\begin{enumerate} \item $(-\cotensor{\coring{C}} X, -
\cotensor{\coring{D}}\Lambda)$ is a Frobenius pair, with
$X_{\coring{D}}$ and $\Lambda_{\coring{C}}$ coflat; \item
$(\Lambda \cotensor{\coring{C}} -, X \cotensor{\coring{D}}-)$ is a
Frobenius pair, with ${}_{\coring{C}}X$ and
${}_{\coring{D}}\Lambda$ coflat; \item $X$ and $\Lambda$ are
coflat quasi-finite injectors on both sides, and $X\simeq
\cohom{\coring{C}}{\Lambda}{\coring{C}}$ in
$^{\coring{C}}\mathcal{M}^{\coring{D}}$ and $\Lambda\simeq
\cohom{\coring{D}}{X}{\coring{D}}$ in
$^{\coring{D}}\mathcal{M}^{\coring{C}}$.
 \end{enumerate}
  If moreover $\coring{C}$ and $\coring{D}$ are coseparable (resp. $A$
and $B$ are von Neumann regular rings), then the following
statements are equivalent
\begin{enumerate}
\item $(-\cotensor{\coring{C}} X,- \cotensor{\coring{D}}\Lambda)$
is a Frobenius pair; \item $(\Lambda \cotensor{\coring{C}} -, X
\cotensor{\coring{D}}-)$ is a Frobenius pair;\item $X$ and
$\Lambda$ are quasi-finite (resp. quasi-finite injectors) on both
sides, and $X\simeq \cohom{\coring{C}}{\Lambda}{\coring{C}}$ in
$^{\coring{C}}\mathcal{M}^{\coring{D}}$ and $\Lambda\simeq
\cohom{\coring{D}}{X}{\coring{D}}$ in
$^{\coring{D}}\mathcal{M}^{\coring{C}}$.
\end{enumerate}
\end{theorem}

\begin{remark}\label{damodulos}
In the case of rings (i.e., $\coring{C} = A$ and $\coring{D} =
B$, Theorem \ref{20a} and the second part of Theorem
\ref{leftrightFrobenius} give \cite[Theorem
2.1]{Castano/Gomez/Nastasescu:1999}. To see this, observe that
${}_AX_B$ is quasi-finite as a right $B$-module if and only if -
$\tensor{A} X : \rmod{A} \rightarrow \rmod{B}$ has a left adjoint,
that is, if and only if ${}_AX$ is finitely generated and
projective. In such a case, the left adjoint is $-
\tensor{B}\hom{A}{X}{A} : \rmod{B} \rightarrow \rmod{A}$. Of
course, $- \cotensor{A} X = - \tensor{A} X$.

The dual characterization in the framework of coalgebras over
fields \cite[Theorem 3.3]{Castano/Gomez/Nastasescu:1999} will be
deduced in Section \ref{condualidad} (see Remarks
\ref{dacomodulos}).
\end{remark}

\section{Frobenius functors between corings with a
duality}\label{condualidad}

We will look to Frobenius functors for corings closer to
coalgebras over fields, in the sense that the categories of
comodules share a fundamental duality.

\par

An object $M$ of a Grothendieck category $\cat{C}$ is said to be
\emph{finitely generated} \cite[p. 121]{Stenstrom:1975} if
whenever $M = \sum_i M_i$ is a direct union of subobjects $M_i$,
then $M = M_{i_0}$ for some index $i_0$. Alternatively, $M$ is
finitely generated if the functor $\hom{\cat{C}}{M}{-}$ preserves
direct unions \cite[Proposition V.3.2]{Stenstrom:1975}. The
category $\cat{C}$ is \emph{locally finitely generated} if it has
a family of finitely generated generators. Recall from \cite[p.
122]{Stenstrom:1975} that a finitely generated object $M$ is
\emph{finitely presented} if every epimorphism $L \rightarrow M$
with $L$ finitely generated has finitely generated kernel. By
\cite[Proposition V.3.4]{Stenstrom:1975}, if $\cat{C}$ is locally
finitely generated, then $M$ is finitely presented if and only if
$\hom{\cat{C}}{M}{-}$ preserves direct limits. For the notion of a
locally projective module we refer to
\cite{Zimmermann-Huisgen:1976}.

\begin{lemma}\label{locally}
Let $\coring{C}$ be a coring over a ring $A$ such that
${}_A\coring{C}$ is flat.
\begin{enumerate}[(1)]
\item\label{fg} A comodule $M \in \rcomod{\coring{C}}$ is finitely
generated if and only if $M_A$ is finitely generated.
\item\label{fp} A comodule $M \in \rcomod{\coring{C}}$ is finitely
presented if $M_A$ is finitely presented. The converse is true
whenever $\rcomod{\coring{C}}$ is locally finitely generated.
\item\label{lfg} If ${}_A\coring{C}$ is locally projective, then
$\rcomod{\coring{C}}$ is locally finitely generated.
\end{enumerate}
\end{lemma}

\begin{proof}
The forgetful functor $U : \rcomod{\coring{C}} \rightarrow
\rmod{A}$  has an exact left adjoint $ - \tensor{A} \coring{C} :
\rmod{A} \rightarrow \rcomod{\coring{C}}$ which preserves direct
limits. Thus, $U$ preserves finitely generated objects and, in
case that $\rcomod{\coring{C}}$ is locally finitely generated,
finitely presented objects. Now, if $M \in \rcomod{\coring{C}}$ is
finitely generated as a right $A$-module, and $M = \sum_i M_i$ as
a direct union of subcomodules, then $U(M) = U(\sum_iM_i) =
\sum_iU(M_i)$, since $U$ is exact and preserves coproducts.
Therefore, $U(M) = U(M_{i_0})$ for some index $i_0$ which implies,
being $U$ a faithfully exact functor, that $M = M_{i_0}$. Thus,
$M$ is a finitely generated comodule. We have thus proved
(\ref{fg}), and the converse to (\ref{fp}). Now, if $M \in
\rcomod{\coring{C}}$ is such that $M_A$ is finitely presented,
then for every exact sequence $0 \rightarrow K \rightarrow L
\rightarrow M \rightarrow 0$ in $\rcomod{\coring{C}}$ with $L$
finitely generated, we get an exact sequence $0 \rightarrow K_A
\rightarrow L_A \rightarrow M_A \rightarrow 0$ with $M_A$ finitely
presented. Thus, $K_A$ is finitely generated and, by (\ref{fg}),
$K \in \rcomod{\coring{C}}$ is finitely generated. This proves
that $M$ is a finitely presented comodule. Finally, (\ref{lfg}) is
a consequence of (\ref{fg}) and
\cite[19.12(1)]{Brzezinski/Wisbauer:2003}.
\end{proof}

The notation $\cat{C}_f$ stands for the full subcategory of a
Grothendieck category $\cat{C}$ whose objects are the finitely
generated objects. The category $\cat{C}$ is locally noetherian
\cite[p. 123]{Stenstrom:1975} if it has a family of noetherian
generators or equivalently, if $\mathbf{C}$ is locally finitely
generated and every finitely generated object of $\mathbf{C}$ is
noetherian. By \cite[Proposition V.4.2, Proposition V.4.1, Lemma
V.3.1(i)]{Stenstrom:1975}, in an arbitrary Grothendieck category,
every finitely generated object is noetherian if and only if every
finitely generated object is finitely presented. The version for
categories of modules of the following result is well-known.

\begin{lemma}\label{10}
Let $\cat{C}$ be a locally finitely generated category.
\begin{enumerate}[(1)]
\item\label{additive} The category $\cat{C}_f$ is additive.
\item\label{cokernels} The category $\cat{C}_f$ has cokernels, and
every monomorphism in $\cat{C}_f$ is a monomorphism in $\cat{C}$.
\item\label{locnoeth} The following statements are equivalent:
  \begin{enumerate}[(a)]
   \item\label{haskernels} The category $\cat{C}_f$ has kernels;
   \item\label{locnoeth2} $\cat{C}$ is locally noetherian;
   \item\label{abelian} $\cat{C}_f$ is abelian;
   \item\label{abeliansub} $\cat{C}_f$ is an abelian subcategory of
$\cat{C}$.
  \end{enumerate}
\end{enumerate}
\end{lemma}

\begin{proof}
(1) Straightforward.

 (2) That $\mathbf{C}_{f}$ has cokernels is
straightforward from \cite[Lemma V.3.1(i)]{Stenstrom:1975}. Now, let
$f:M\rightarrow N$ be a monomorphism in $\mathbf{C}_{f}$ and
$\xi:X\rightarrow M$ be a morphism in $\mathbf{C}$ such that
$f\xi=0.$ Suppose that $X=\bigcup_{i\in I}X_i$, where
$X_i\in\mathbf{C}_{f}$, and $\iota_i:X_i\rightarrow X$, $i\in I$ the
canonical injections$.$ Then $f\xi\iota_i=0$, and $\xi\iota_i=0$,
for every $i$, and by the definition of the inductive limit, $\xi=0.$%

(3) $(b)\Rightarrow(a)$ Straightforward from \cite[Proposition
V.4.1]{Stenstrom:1975}.

 $(d)\Rightarrow(c)$ and $(c)\Rightarrow(a)$
are trivial.

$(a)\Rightarrow(b)$ Let $M\in\mathbf{C}_{f}$, and $K$ be a subobject
of $M.$ Let $\iota :L\rightarrow M$ the kernel of the canonical
morphism $f:M\rightarrow M/K$ in $\mathbf{C}_{f}.$ Suppose that
$K=\bigcup_{i\in I}K_i$, where $K_i\in\mathbf{C}_{f}$, for every
$i\in I.$ By the universal property of the kernel, there exist a
unique morphism $\alpha:L\rightarrow K,$ and a unique morphism
$\beta_i:K_i\rightarrow L$, for every $i\in I$, making commutative
the diagrams
\[%
\xymatrix{
&   K_i \ar[d] \ar[dl]_{\beta_i} & \\
L \ar[dr]_{\alpha} \ar[r]^{\iota} & M \ar[r]^{f} & M/K \\
   & K \ar[u] & .}
\]
By (2), $\iota$ is a monomorphism in $\mathbf{C}$, then for every
$K_i\subset K_{j}$, the diagram
\[%
\xymatrix{
K_i  \ar[d] \ar[r]^{\beta_i} & L\\
\ar[ur]_{\beta_{j}}  K_{j} &   }
\]
commutes. Therefore we have the commutative diagram
\[%
\xymatrix{
& K \ar[d] \ar[dl]_{{\lim \atop {\longrightarrow \atop I}}\beta_i}\\
L \ar[r]^{\iota} & M.}
\]
Then $K\simeq L$, and hence $K\in\mathbf{C}_{f}.$ Finally, by
\cite[Proposition V.4.1]{Stenstrom:1975}, $M$ is noetherian in
$\mathbf{C.}$

$(b)\Rightarrow(d)$ Straightforward from \cite[Theorem
3.41]{Freyd:1964}.
\end{proof}

The following generalization of \cite[Proposition
3.1]{Castano/Gomez/Nastasescu:1999} will allow us to give an
alternative proof to the equivalence ``$(1)\Leftrightarrow(4)$'' of
Theorem \ref{leftrightFrobeniusduality} (see Remarks
\ref{dacomodulos}).

\begin{proposition}\label{11}
Let $\mathbf{C}$ and $\mathbf{D}$ be two locally noetherian
categories. Then
\begin{enumerate}[(1)]
\item If $F:\mathbf{C}\rightarrow\mathbf{D}$ is a Frobenius
functor, then its restriction $F_f:\mathbf{C}_{f}\rightarrow
\mathbf{D}_{f}$ is a Frobenius functor.
\item If $H:\mathbf{C}%
_{f}\rightarrow\mathbf{D}_{f}$ is a Frobenius functor, then $H$
can be
uniquely extended to a Frobenius functor $\overline{H}:\mathbf{C}%
\rightarrow\mathbf{D}.$ \item The assignment $F\mapsto F_f$ defines
a bijective correspondence (up to natural isomorphisms) between
Frobenius functors from $\mathbf{C}$ to $\mathbf{D}$ and Frobenius
functors from $\mathbf{C}_{f}$ to $\mathbf{D}_f.$
\item In particular, if
$\mathbf{C}=\mathcal{M}^{\coring{C}}$ and $\mathbf{D}=\mathcal{M}%
^\coring{D}$ are locally noetherian such that $_A\coring{C}$ and
$_B\coring{D}$ are flat, then $F:\mathcal{M}^{\coring{C}%
}\rightarrow\mathcal{M}^\coring{D}$ is a Frobenius functor if and
only if it preserves direct limits and comodules which are finitely
generated as right $A$-modules,
and the restriction functor $F_f:\mathcal{M}_{f}^{\coring{C}}%
\rightarrow\mathcal{M}_{f}^{\coring{D}}$ is a Frobenius functor.
\end{enumerate}
\end{proposition}

\begin{proof}
The proofs of \cite[Proposition 3.1 and Remark
3.2]{Castano/Gomez/Nastasescu:1999} remain valid for our
situation, but with some minor modifications: to prove that
$\overline{H}$ is well-defined, we use Lemma \ref{10}. In the
proof of the statements (1), (2) and (3) we use the Grothendieck
AB 5 condition.
\end{proof}

In order to generalize \cite[Proposition A.2.1]{Takeuchi:1977ca}
and its proof, we need the following lemma.

\begin{lemma}\label{12}
\begin{enumerate}[(1)] \item Let $\mathbf{C}$ be a locally noetherian 
category, let
$\mathbf{D}$ be an arbitrary Grothendieck category,
$F:\mathbf{C}\rightarrow\mathbf{D}$ be an
arbitrary functor which preserves direct limits, and $F_f:\mathbf{C}%
_{f}\rightarrow\mathbf{D}$ be its restriction to $\mathbf{C}_{f}.$
Then $F$ is exact (faithfully exact, resp. left, right exact) if and
only if $F_f$ is exact (faithfully exact, resp. left, right exact).

 In particular, an object $M$ in $\cat{C}_f$ is projective
(resp. projective generator) if and only if it is projective
(resp. projective generator) in $\cat{C}$.

\item Let $\mathbf{C}$ be a locally noetherian category. For every
object $M$ of $\mathbf{C}$, the following conditions are equivalent
\begin{enumerate}[(a)]
\item $M$ is injective (resp. an injective cogenerator); \item the
contravariant functor $\hom{\mathbf{C}}{-}{M}
:\mathbf{C}\rightarrow \mathbf{Ab}$ is exact (resp. faithfully
exact); \item the contravariant
functor $\hom{\mathbf{C}}{-}{M}_{f}:\mathbf{C}_{f}%
\rightarrow\mathbf{Ab}$ is exact (resp. faithfully exact).
\end{enumerate}
 In particular, an object $M$ in $\cat{C}_f$ is injective (resp.
injective cogenerator) if and only if it is injective (resp.
injective cogenerator) in $\cat{C}$.
\end{enumerate}
\end{lemma}

\begin{proof}
(1) The ``only if'' part is straightforward from the fact that the
injection functor $\mathbf{C}_f\rightarrow\mathbf{C}$ is faithfully
exact.

 For the ``if'' part, suppose
that $F_f$ is left exact. Let $f:M\rightarrow N$ be a morphism in
$\mathbf{C}$. Put $M=\bigcup_{i\in I}M_i$ and $N=\bigcup_{j\in
J}N_{j},$ as direct union of directed families of
finitely generated subobjects. For $(i,j)\in I\times J,$ let $M_{i,j}%
=M_i \cap f^{-1}(N_j)$, and $f_{i,j}:M_{i,j}\rightarrow N_{j}$ be
the restriction of $f$ to $M_{i,j}.$ We have
$f=\underset{\underset{I\times J}{\longrightarrow}}{\lim}f_{i,j}$
and then $F(f)=\underset {\underset{I\times
J}{\longrightarrow}}{\lim}F_f(f_{i,j}).$ Hence

 $\ker F(f)=\ker\underset{\underset{I\times
J}{\longrightarrow}}{\lim}F_f(f_{i,j})=\underset {\underset{I\times
J}{\longrightarrow}}{\lim}\ker F_f( f_{i,j})
=\underset{\underset{I\times J}{\longrightarrow}}{\lim}F_f(\ker
f_{i,j})=\underset{\underset{I\times
J}{\longrightarrow}}{\lim}F(\ker f_{i,j})$ (by Lemma \ref{10})
$=F(\underset{\underset{I\times
J}%
{\longrightarrow}}{\lim}\ker f_{i,j})=F(\ker f).$

 Finally $F$ is left exact. Analogously, it can be proved that $F_f$ is 
right exact
implies that $F$ is also right exact. Now, suppose that $F_f$ is
faithfully exact. We have already proved that $F$ is exact. It
remains to prove that $F$ is faithful. For this, let $0\neq
M=\bigcup_{i\in I}M_i$ be an object of $\mathbf{C}$, where $M_i$ is
finitely generated for every $i\in I.$ We have
\[
F(M)=\underset{\underset{I}{\longrightarrow}}{\lim}F_f(M_i)
\simeq\sum_iF_f(M_i)
\]
(since $F$ is exact). Since $M\neq0$, there exists some $i_{0}\in I$
such that $M_{i_{0}}\neq0.$ By \cite[Proposition
IV.6.1]{Stenstrom:1975}, $F_f(M_{i_{0}})\neq0,$ hence $F(M)\neq0.$
Also by \cite[Proposition IV.6.1]{Stenstrom:1975}, $F$ is faithful.

(2) $(a)\Leftrightarrow(b)$ Obvious.

 $(b)\Rightarrow(c)$
Analogous to that of the ``only if'' part of (1).

$(c)\Rightarrow(a)$ That $M$ is injective is a consequence of
\cite[Proposition V.2.9, Proposition V.4.1]{Stenstrom:1975}. Now,
suppose moreover that $\hom{\mathbf{C}}{-}{M}_{f}$ is faithful. Let
$L$ be a non-zero object of $\mathbf{C}$, and $K$ be a non-zero
finitely generated subobject of $L.$ By \cite[Proposition
IV.6.1]{Stenstrom:1975}, there exists a non-zero morphism
$K\rightarrow M$. Since $M$ is injective, there exists a non-zero
morphism $L\rightarrow M$ making commutative the following diagram
\[%
\xymatrix{
&   M   & \\
0  \ar[r] & K \ar[r] \ar[u] & L. \ar[ul]}
\]
From \cite[Proposition IV.6.5]{Stenstrom:1975}, it follows that
$M$\ is a cogenerator.
\end{proof}

If $\coring{C}_A$ is flat and $M \in \rcomod{\coring{C}}$ is
finitely presented as right $A$-module, then
\cite[19.19]{Brzezinski/Wisbauer:2003} the dual left $A$-module
$\rdual{M} = \hom{A}{M}{A}$ has a left $\coring{C}$-comodule
structure
$$\rdual{M} \simeq
\hom{\coring{C}}{M}{\coring{C}} \subseteq \hom{A}{M}{\coring{C}}
\simeq \coring{C} \tensor{A} \rdual{M}.$$ Now, if ${}_A\rdual{M}$
turns out to be finitely presented and ${}_A\coring{C}$ is flat,
then $\ldual{(\rdual{M})} = \hom{A}{\rdual{M}}{A}$ is a right
$\coring{C}$-comodule and the canonical map $\sigma_M : M
\rightarrow \ldual{(\rdual{M})}$ is a homomorphism in
$\rcomod{\coring{C}}$. This construction leads to a duality (i.e.
a contravariant equivalence)
\[
\rDual: \rcomod{\coring{C}}_0 \leftrightarrows
\lcomod{\coring{C}}_0: \lDual
\]
between the full subcategories $\rcomod{\coring{C}}_0$ and
$\lcomod{\coring{C}}_0$ of $\rcomod{\coring{C}}$ and
$\lcomod{\coring{C}}$ whose objects are the comodules which are
finitely generated and projective over $A$ on the corresponding
side (this holds even without flatness assumptions of
$\coring{C}$). Call it the \emph{basic duality} (details may be
found in \cite{Caenepeel/DeGroot/Vercruysse:unp}). Of course, in
the case that $A$ is semisimple (e.g. for coalgebras over fields)
these categories are that of finitely generated comodules, and
this basic duality plays a remarkable role in the study of several
notions in the coalgebra setting (e.g. Morita equivalence
\cite{Takeuchi:1977}, semiperfect coalgebras \cite{Lin:1977},
Morita duality \cite{Gomez/Nastasescu:1995},
\cite{Gomez/Nastasescu:1996}, or Frobenius Functors
\cite{Castano/Gomez/Nastasescu:1999}). It would be interesting to
know, in the coring setting, to what extent the basic duality can
be extended to the subcategories $\rcomod{\coring{C}}_f$ and
$\lcomod{\coring{C}}_f$, since, as we will try to show in this
section, this allows to obtain better results. Of course, this is
the underlying idea when the ground ring $A$ is assumed to be
Quasi-Frobenius (see \cite{ElKaoutit/Gomez:2002unp} for the case
of semiperfect corings and Morita duality), but we hope future
developments of the theory will be aided by the more general
setting we propose here.

Consider contravariant functors between Grothendieck categories $
H: \cat{A} \leftrightarrows \cat{A}': H', $ together with natural
transformations $\tau : 1_{\cat{A}} \rightarrow H' \circ H$ and
$\tau' : 1_{\cat{A}'} \rightarrow H \circ H'$, satisfying the
condition $H(\tau_X) \circ \tau'_{H(X)} = 1_{H(X)}$ and
$H'(\tau'_{X'}) \circ \tau_{H'(X')} = 1_{H'(X')}$ for $X \in
\cat{A}$ and $X' \in \cat{A}'$. Following
\cite{Colby/Fuller:1983}, this situation is called a \emph{right
adjoint pair}.

\begin{proposition}\label{lnoethrightadj}
Let $\coring{C}$ be an $A$-coring such that ${}_A\coring{C}$ and
$\coring{C}_A$ are flat. Assume that $\rcomod{\coring{C}}$ and
$\lcomod{\coring{C}}$ are locally noetherian categories. If
${}_A\rdual{M}$ and $\ldual{N}_A$ are finitely generated modules
for every $M \in \rcomod{\coring{C}}_f$ and $N \in
\lcomod{\coring{C}}_f$, then the basic duality extends to a right
adjoint pair $ \rDual: \rcomod{\coring{C}}_f \leftrightarrows
\lcomod{\coring{C}}_f: \lDual$.
\end{proposition}

\begin{proof}
If $M \in \rcomod{\coring{C}}_f$ then, since $\rcomod{\coring{C}}$
is locally noetherian, $M_{\coring{C}}$ is finitely presented. By
Lemma \ref{locally}, $M_A$ is finitely presented and the left
$\coring{C}$-comodule $\rdual{M}$ makes sense. Now, the assumption
${}_A\rdual{M}$ finitely generated implies, by Lemma
\ref{locally}, that $\rdual{M} \in \lcomod{\coring{C}}_f$. We have
then the functor $\rDual : \rcomod{\coring{C}}_f \rightarrow
\lcomod{\coring{C}}_f$. The functor $\rDual$ is analogously
defined, and the rest of the proof consists of straightforward
verifications.
\end{proof}

\begin{example}
The hypotheses are fulfilled if ${}_A\coring{C}$ and
$\coring{C}_A$ are locally projective and $A$ is left and right
noetherian (in this case the right adjoint pair already appears in
\cite{ElKaoutit/Gomez:2002unp}). But there are situations in which
no finiteness condition need to be required to $A$: this is the
case, for instance, of cosemisimple corings (see \cite[Theorem
3.1]{ElKaoutit/Gomez/Lobillo:2001unp}). In particular, if an
arbitrary ring $A$ contains a division ring $B$, then, by
\cite[Theorem 3.1]{ElKaoutit/Gomez/Lobillo:2001unp} the canonical
coring $A \tensor{B} A$ satisfies all hypotheses in Proposition
\ref{lnoethrightadj}.
\end{example}

\begin{definition}
Let $\coring{C}$ be a coring over $A$ satisfying the assumptions
of Proposition \ref{lnoethrightadj}. We will say that $\coring{C}$
\emph{has a duality} if the basic duality extends to a duality $$
\rDual: \rcomod{\coring{C}}_f \leftrightarrows
\lcomod{\coring{C}}_f: \lDual.$$
\end{definition}

We have the following examples of a coring which has a
duality:\begin{itemize} \item $\coring{C}$ is a coring over a QF
ring $A$ such that ${}_A\coring{C}$ and $\coring{C}_A$ are flat (and
hence projective);\item $\coring{C}$ is a cosemisimple coring;
where, by \cite[Theorem 3.1]{ElKaoutit/Gomez/Lobillo:2001unp},
$\rcomod{\coring{C}}_f$ and $\lcomod{\coring{C}}_f$ are equal to
$\rcomod{\coring{C}}_0$ and $\lcomod{\coring{C}}_0$, respectively;
\item $\coring{C}$ is a coring over $A$ such that ${}_A\coring{C}$
and $\coring{C}_A$ are flat and semisimple, $\rcomod{\coring{C}}$
and $\lcomod{\coring{C}}$ are locally noetherian categories, and the
dual of every simple right (resp. left) $A$-module in the
decomposition of $\coring{C}_A$ (resp. ${}_A\coring{C}$) as a direct
sum of simple $A$-modules is finitely generated and $A$-reflexive
(in fact, every right (resp. left) $\coring{C}$-comodule $M$ becomes
a submodule of the semisimple right (resp. left) $A$-module
$M\otimes_A\coring{C}$, and hence $M_A$ (resp. ${}_AM$) is also
semisimple).\end{itemize}

\begin{proposition}\label{13}
Suppose that the coring $\coring{C}$ has a duality. Let
$M\in\mathcal{M}^{\coring{C}}$ such that $M_A$ is flat. The
following are equivalent
\begin{enumerate} \item $M$ is coflat (resp. faithfully coflat);
\item
$\hom{\coring{C}}{-}{M}_{f}:\mathcal{M}_{f}^{\coring{C}}\rightarrow\mathcal{M}_{k}$
is exact (resp. faithfully exact); \item $M$ is injective (resp. an
injective cogenerator).\end{enumerate}
\end{proposition}

\begin{proof}
Let $M\in\mathcal{M}%
^{\coring{C}}$ and $N\in{}^{\coring{C}}\mathcal{M}_{f}.$ We have
the following commutative diagram (in $\mathcal{M}_{k}$)
\[%
\xymatrix{  0 \ar[r] & M\square_{\coring{C}}N  \ar[r] &
M\otimes_AN \ar[rrr]^{\rho_M\otimes_AN-M\otimes_A\lambda_N%
} \ar[d] & & & M\otimes_A\coring{C}\otimes_AN \ar[d] \\
 0 \ar[r] & \hom{\coring{C}}{N^*}{M}\ar[r] &
\hom{A}{N^*}{M}\ar[rrr]^{f\mapsto\rho_Mf-(f\otimes_A\coring{C})\rho_{N^*}
} & & & \hom{A}{N^*}{M\otimes_A\coring{C}},}
\]
where the vertical maps are the canonical maps. By the universal
property of the kernel, there is a unique morphism
$\eta_{M,N}:M\square_{\coring{C}}N\rightarrow
\hom{\coring{C}}{N^*}{M}$ making commutative the above diagram. By
the cube Lemma (see \cite[Proposition II.1.1]{Mitchell:1965}),
$\eta$ is a natural transformation of bifunctors. If $M_A$ is flat
then $\eta_{M,N}$ is an isomorphism for every
$N\in{}^{\coring{C}}\mathcal{M}_{f}.$ We have
\[
M\square_{\coring{C}}-\simeq \hom{\coring{C}}{-}{M}_{f} \circ
(-)^*:{}^{\coring{C}}\mathcal{M}_{f}\rightarrow \mathcal{M}_{k}.
\]
Then, by Lemma \ref{12}, $M_{\coring{C}}$ is coflat (resp.
faithfully coflat) iff
$M\square_{\coring{C}}-:{}^{\coring{C}}\mathcal{M}_{f}\rightarrow
\mathcal{M}_{k}$ is exact (resp. faithfully exact)  iff
$\hom{\coring{C}}{-}{M}_{f}:\mathcal{M}_{f}^{\coring{C}%
}\rightarrow\mathcal{M}_{k}$ is exact (resp. faithfully exact) iff
$M_{\coring{C}}$ is injective (resp. an injective cogenerator).
\end{proof}

The particular case of the following result for coalgebras over
commutative ring is given in \cite{AlTakhman:2002}.

\begin{corollary}\label{injectiveinjector}
Let $N\in{}^{\coring{C}}\mathcal{M}^{\coring{D}}$ be a bicomodule.
Suppose that $A$ is a QF ring and $\coring{D}$ has a duality. If
$N$ is injective in $\mathcal{M}%
^{\coring{D}}$ such that $N_B$ is flat, then $N$ is an injector as a
right $\coring{D}$-comodule.
\end{corollary}

\begin{proof}
Let $X_A$ be an injective module. Since $A$ is a QF ring, $X_A$ is
projective. We have then the natural isomorphism
\[
(X\otimes_AN)\square_{\coring{D}}-\simeq
X\otimes_A(N\square_{\coring{D}}-)
:{}^{\coring{D}}\mathcal{M}\rightarrow \mathcal{M}_{k}.
\]
By Proposition \ref{13}, $N_{\coring{D}}$ and $X_A$ are coflat, and
then $X\otimes_AN$ is coflat. Now, since $X\otimes_AN$ is a flat
right $B$-module, and by Proposition \ref{13}, $X\otimes_AN$ is
injective in $\mathcal{M}^{\coring{D}}$.
\end{proof}

 The last two results allow to improve our general statements in
 Section \ref{Frobeniusgeneral} for corings having a duality.

\begin{proposition}\label{adjointpairduality}
Suppose that $\coring{C}$ and $\coring{D}$ have a duality.
Consider the following statements
\begin{enumerate} \item $(-\cotensor{\coring{C}} X, -
\cotensor{\coring{D}}\Lambda)$ is an adjoint pair of functors, with
$_AX$ and $\Lambda_A$ flat;\item $\Lambda$ is quasi-finite injective
as a right ${\coring{C}}$-comodule, with $_AX$ and $\Lambda_A$ flat
and $X\simeq \cohom{\coring{C}}{\Lambda}{\coring{C}}$ in
$^{\coring{C}}\mathcal{M}^{\coring{D}}$.
\end{enumerate}
We have (1) implies (2), and the converse is true if in particular
$B$ is a QF ring.
\end{proposition}

\begin{proof}
$(1)\Rightarrow(2)$ From Proposition \ref{13}, $\coring{D}$ is
injective in $\mathcal{M}^{\coring{D}}$. Since the functor
$-\cotensor{\coring{C}} X$ is exact, $\Lambda\simeq
\coring{D}\cotensor{\coring{D}} \Lambda$ is injective in
$\mathcal{M}^{\coring{C}}$.

 $(2)\Rightarrow(1)$ Assume
that $B$ is a QF ring. From Corollary \ref{injectiveinjector},
$\Lambda$ is quasi-finite injector as a right
${\coring{C}}$-comodule, and Proposition \ref{19} achieves the
proof.
\end{proof}

We are now in a position to state and prove our main result of this
section.

\begin{theorem}\label{leftrightFrobeniusduality}
Suppose that $\coring{C}$ and $\coring{D}$ have a duality. Let
$X\in{}^{\coring{C}}\mathcal{M}%
^{\coring{D}}$ and
$\Lambda\in{}^{\coring{D}}\mathcal{M}^{\coring{C}}.$ The following
statements are equivalent \begin{enumerate} \item
$(-\cotensor{\coring{C}} X, - \cotensor{\coring{D}}\Lambda)$ is a
Frobenius pair, with $X_B$ and $\Lambda_A$ flat; \item $(\Lambda
\cotensor{\coring{C}} -, X \cotensor{\coring{D}}-)$ is a Frobenius
pair, with $_AX$ and $_B\Lambda$ flat; \item $X$ and $\Lambda$ are
quasi-finite injector on both sides, and $X\simeq
\cohom{\coring{C}}{\Lambda}{\coring{C}}$ in
$^{\coring{C}}\mathcal{M}^{\coring{D}}$ and $\Lambda\simeq
\cohom{\coring{D}}{X}{\coring{D}}$ in
$^{\coring{D}}\mathcal{M}^{\coring{C}}$.

In particular, if $A$ and $B$ are QF rings, then the above
statements are equivalent to \item $X$ and $\Lambda$ are
quasi-finite injective on both sides, and $X\simeq
\cohom{\coring{C}}{\Lambda}{\coring{C}}$ in
$^{\coring{C}}\mathcal{M}^{\coring{D}}$ and $\Lambda\simeq
\cohom{\coring{D}}{X}{\coring{D}}$ in
$^{\coring{D}}\mathcal{M}^{\coring{C}}$.
\end{enumerate}
Finally, suppose that $\coring{C}$ and $\coring{D}$ are
cosemisimple corings. Let
$X\in{}^{\coring{C}}\mathcal{M}%
^{\coring{D}}$ and
$\Lambda\in{}^{\coring{D}}\mathcal{M}^{\coring{C}}.$ The following
statements are equivalent \begin{enumerate} \item
$(-\cotensor{\coring{C}} X, - \cotensor{\coring{D}}\Lambda)$ is a
Frobenius pair; \item $(\Lambda \cotensor{\coring{C}} -, X
\cotensor{\coring{D}}-)$ is a Frobenius pair; \item $X$ and
$\Lambda$ are quasi-finite on both sides, and $X\simeq
\cohom{\coring{C}}{\Lambda}{\coring{C}}$ in
$^{\coring{C}}\mathcal{M}^{\coring{D}}$ and $\Lambda\simeq
\cohom{\coring{D}}{X}{\coring{D}}$ in
$^{\coring{D}}\mathcal{M}^{\coring{C}}$.
\end{enumerate}
\end{theorem}

\begin {proof}
We start by proving the first part. In view of Theorem
\ref{leftrightFrobenius} and Theorem \ref{20a} it suffices to show
that if $(-\cotensor{\coring{C}} X, - \cotensor{\coring{D}}\Lambda)$
is a Frobenius pair, the condition ``$X_{\coring{D}}$ and
$\Lambda_{\coring{C}}$ are coflat'' is equivalent to ``$X_B$ and
$\Lambda_A$  are flat''. Indeed, the first implication is obvious,
for the converse, assume that $X_B$ and $\Lambda_A$ are flat. By
Proposition \ref{adjointpairduality}, $X$ and $\Lambda$ are
injective in $\mathcal{M}^{\coring{D}}$ and
$\mathcal{M}^{\coring{C}}$ respectively, and they are coflat by
Proposition \ref{13}. The particular case is straightforward from
Proposition \ref{adjointpairduality} and the above equivalences .
\\ Now we will show the second part. We know that cosemisimple corings 
have a duality. By \cite[Theorem
3.1]{ElKaoutit/Gomez/Lobillo:2001unp}, every comodule category
over a cosemisimple coring is a spectral category (see \cite[p.
128]{Stenstrom:1975}). Thus, the bicomodules
${}_{\coring{C}}X_{\coring{D}}$ and
${}_{\coring{D}}\Lambda_{\coring{C}}$ are coflat and injector on
both sides (we can see this directly by using the fact that every
additive functor between abelian categories preserves split
exactness). Now, apply the first part.
\end{proof}

\begin{remarks}\label{dacomodulos}
\begin{enumerate}
\item The equivalence ``$(1)\Leftrightarrow(4)$'' of the last
theorem is a generalization of \cite[Theorem
3.3]{Castano/Gomez/Nastasescu:1999}. The proof of \cite[Theorem
3.3]{Castano/Gomez/Nastasescu:1999} gives an alternative proof of
``$(1)\Leftrightarrow(4)$'' of Theorem
\ref{leftrightFrobeniusduality}, using Proposition \ref{11}. \item
The adjunction of Proposition \ref{19a} and Proposition
\ref{adjointpairduality} generalizes the coalgebra version of
Morita's theorem \cite[Theorem
4.2]{Caenepeel/DeGroot/Militaru:2002}.
\end{enumerate}
\end{remarks}

\begin{example}
Let $A$ be a $k$-algebra. Put $\coring{C}=A$ and $\coring{D}=k.$
The bicomodule $A\in{}^{\coring{C}}\mathcal{M}^{\coring{D}}$ is
quasi-finite as a right $\coring{D}$-comodule. $A$ is an injector
as a right $\coring{D} $-comodule if and only if the $k$-module
$A$ is flat. If we take $A=k=\mathbb{Z} $, the bicomodule $A$ is
quasi-finite and injector as a right $\coring{D}$-comodule but it
is not injective in $\mathcal{M}^{\coring{D}}.$ Hence, the
assertion ``$(-\cotensor{\coring{C}} X, -
\cotensor{\coring{D}}\Lambda)$ is an adjoint pair of functors''
does not imply in general the assertion ``$\Lambda$ is
quasi-finite injective as a right ${\coring{C}}$-comodule and
$X\simeq \cohom{\coring{C}}{\Lambda}{\coring{C}}$ in
$^{\coring{C}}\mathcal{M}^{\coring{D}}$'', and the following
statements are not equivalent in general:
\begin{enumerate} \item $(-\cotensor{\coring{C}} X, -
\cotensor{\coring{D}}\Lambda)$ is a Frobenius pair; \item $X$ and
$\Lambda$ are quasi-finite injective on both sides, and $X\simeq
\cohom{\coring{C}}{\Lambda}{\coring{C}}$ in
$^{\coring{C}}\mathcal{M}^{\coring{D}}$ and $\Lambda\simeq
\cohom{\coring{D}}{X}{\coring{D}}$ in
$^{\coring{D}}\mathcal{M}^{\coring{C}}$.
\end{enumerate}
On the other hand, there exists a commutative self-injective ring
wich is not coherent. By a theorem of S.U.~Chase (see for example
\cite[Theorem 19.20]{Anderson/Fuller:1992}), there exists then a
$k$-algebra $A$ which is injective, but not flat as $k$-module.
Hence, the bicomodule
$A\in{}^{\coring{C}}\mathcal{M}^{\coring{D}}$ is quasi-finite and
injective as a right $\coring{D}$-comodule, but not an injector as
a right $\coring{D}$-comodule.
\end{example}

\section{Applications to induction functors}\label{Frobmor}

We start this section by recalling from \cite{Gomez:2002}, that a
coring homomorphism from the coring $\coring{C}$\ into the coring
$\coring{D}$ is a pair $(\varphi,\rho)$, where $\rho:A\rightarrow B$
is a homomorphism of $k$-algebras and
$\varphi:\coring{C}\rightarrow\coring{D}$ is a homomorphism of
$A$-bimodules such that
\[
\epsilon_{\coring{D}}\circ\varphi=\rho\circ\epsilon_{\coring{C}}\qquad\textrm{and
}\qquad\Delta_{\coring{D}}\circ\varphi=\omega_{\coring{D},\coring{D}}\circ(
\varphi\otimes_A\varphi)\circ\Delta_{\coring{C}},
\]
where $\omega_{\coring{D},\coring{D}}:\coring{D}\otimes_A\coring{D}%
\rightarrow\coring{D}\otimes_B\coring{D}$ is the canonical map
induced by $\rho:A\rightarrow B.$

\medskip

Now we will characterize when the induction functor $ - \tensor{A}
B : \rcomod{\coring{C}} \rightarrow \rcomod{\coring{D}}$ defined
in \cite[Proposition 5.3]{Gomez:2002} is a Frobenius functor. The
coaction of $\coring{D}$ over $M \tensor{A} B$ is given, when
expressed in Sweedler's sigma notation, by
\[
\rho_{M \tensor{A} B}(m \tensor{A}b) = \sum m_{(0)} \tensor{A} 1_B
\tensor{B} \varphi(m_{(1)})b,
\]
where $M$ is a right $\coring{C}$-comodule with coaction
$\rho_M(m) = \sum m_{(0)} \tensor{A} m_{(1)}$. We also define the
functor $-\square_{\coring{D}}(B\tensor{A}\coring{C}):
\rcomod{\coring{D}}\rightarrow \rcomod{\coring{C}}$, where the
left coaction on the left $B$-module $B\tensor{A}\coring{C}$ is
given by:
$$\lambda_{B\tensor{A}\coring{C}}:B\tensor{A}\coring{C}\rightarrow
\coring{D}\tensor{B}B\tensor{A}\coring{C}\simeq
\coring{D}\tensor{A}\coring{C},\quad b\tensor{A}c\mapsto \sum
b\varphi(c_{(1)})\tensor{A}c_{(2)},$$ where
$\Delta_{\coring{C}}(c)=\sum c_{(1)}\tensor{A}c_{(2)}$. Moreover,
if $\omega_{Y,B\tensor{A}\coring{C}}$ is ${}_A\coring{C}$--pure
for every right $\coring{D}$--comodule $Y$, then, by
\cite[Proposition 5.4]{Gomez:2002}, we have the adjoint pair of
functors
$(-\otimes_AB,-\square_{\coring{D}}(B\tensor{A}\coring{C}))$.

\begin{theorem}\label{23}
Let $(\varphi,\rho) :\coring{C}\rightarrow\coring{D}$ be a
homomorphism of corings such that $_A\coring{C}$ and $_B\coring{D} $
are flat. The following statements are
equivalent\begin{enumerate}[(a)] \item
$-\otimes_AB:\mathcal{M}^{\coring{C}}\rightarrow\mathcal{M}%
^{\coring{D}}$ is a Frobenius functor;\item the $\coring{C}%
-\coring{D}$-bicomodule $\coring{C}\otimes_AB$ is quasi-finite and
injector as a right $\coring{D}$-comodule and there exists an
isomorphism of $\coring{D}-\coring{C}$-bicomodules
$\cohom{\coring{D}}{\coring{C}\otimes_AB}{\coring{D}}\simeq B\otimes
_A\coring{C}.$
\end{enumerate}
Moreover, if $\coring{C}$ and $\coring{D}$ are coseparable, then
the condition ``injector'' in $(b)$ can be deleted.
\end{theorem}

\begin{proof}
First observe that $-\otimes_AB$ is a Frobenius functor if and only
if $(-\otimes_AB,-\square_{\coring{D}}(B\otimes _A\coring{C}))$ is a
Frobenius pair (by \cite[Proposition 5.4]{Gomez:2002}). A
straightforward computation shows that the map $\rho_M \tensor{A} B
: M \tensor{A} B \rightarrow M \tensor{A} \coring{C} \tensor{A} B$
is a homomorphism of $\coring{D}$-comodules. We have thus a
commutative diagram in $\rcomod{\coring{C}}$ with exact row
\[
\xymatrix{ 0 \ar[r] & M \cotensor{\coring{C}} (\coring{C}
\tensor{A} B) \ar^{\iota}[r] & M \tensor{A} \coring{C} \tensor{A}
B \ar^-{\omega_{M,\coring{C} \tensor{A} B}}[rr] & & M \tensor{A}
\coring{C} \tensor{A} \coring{C} \tensor{A} B \\
 & M \tensor{A} B \ar^{\psi_M}[u] \ar_-{\rho_M \tensor{A} B}[ur] & &
 &},
\]
where $\psi_M$ is defined by the universal property of the kernel.
Since ${}_B\coring{D}$ is flat, to prove that $\psi_M$ is an
isomorphism of $\coring{D}$-comodules it is enough to check that it
is bijective, as the forgetful functor $U : \rcomod{\coring{D}}
\rightarrow \rmod{B}$ is faithfully exact. Some easy computations
show that the map $(M \tensor{A} \epsilon_{\coring{C}} \tensor{A}
B)\circ \iota$ is the inverse in $\rmod{B}$ to $\psi_M$. From this,
we deduce a natural isomorphism $ \psi : -\otimes_AB\simeq
-\square_{\coring{C}}( \coring{C}\otimes_AB)$. The equivalence
between (a) and (b) is then obvious from Proposition \ref{19} and
Proposition \ref{19a}.
\end{proof}

When applied to the case where $\coring{C} = A$ and $\coring{D} =
B$ are the trivial corings (which are separable), Theorem \ref{23}
gives functorial Morita's characterization of Frobenius ring
extensions given in \cite[Theorem 5.1]{Morita:1965}. This follows
from \cite[Example 4.3]{Gomez:2002}: In the case $\coring{C} = A$,
$\coring{D} = B$ we have that $A\tensor{A}B \cong B$ is
quasi-finite as a right $B$--comodule if and only if ${}_AB$ is
finitely generated an projective, and, in this case,
$\cohom{B}{B}{-} \simeq - \tensor{B}\hom{A}{{}_AB}{A}$.

Theorem \ref{24} generalizes the characterization of Frobenius
extension of coalgebras over fields \cite[Theorem
3.5]{Castano/Gomez/Nastasescu:1999}. It is then reasonable to give
the following definition.

\begin{definition}\label{Frobenius extension}
Let $(\varphi,\rho) :\coring{C}\rightarrow\coring{D}$ be a
homomorphism of corings such that $_A\coring{C}$ and $_B\coring{D} $
are flat. It is said to be a \emph{right Frobenius} morphism of
corings if $ - \tensor{A} B : \rcomod{\coring{C}} \rightarrow
\rcomod{\coring{D}}$ is a Frobenius functor.
\end{definition}

\begin{theorem}\label{24}
Suppose that the algebras $A$ and $B$ are QF rings.\newline Let
$(\varphi,\rho):\coring{C}\rightarrow\coring{D}$ be a homomorphism
of corings such that the modules $_A\coring{C}$, $_B\coring{D}$ and
$\coring{D}_B$ are projective. Then the following statements are
equivalent
\begin{enumerate}[(a)]
\item $-\otimes_AB:\mathcal{M}^{\coring{C}}\rightarrow
\mathcal{M}^{\coring{D}}$ is a Frobenius functor;\item the
$\coring{C}-\coring{D}$-bicomodule $\coring{C}\otimes_AB$ is
quasi-finite as a right $\coring{D}$-comodule,
$(\coring{C}\otimes_AB) _{\coring{D}}$ is injective and there exists
an isomorphism of
$\coring{D}%
-\coring{C}$-bicomodules
$\cohom{\coring{D}}{\coring{C}\otimes_AB}{\coring{D}}\simeq
B\otimes_A\coring{C}.$
\end{enumerate}
\end{theorem}

\begin{proof}
Obvious from Proposition \ref{adjointpairduality}.
\end{proof}

Now, suppose that the forgetful functor $\mathcal{M}^{\coring{C}%
}\rightarrow\mathcal{M}_A$ is a Frobenius functor. Then the functor
$-\otimes_A\coring{C}:\mathcal{M}_A\rightarrow\mathcal{M}_A$ is also
a Frobenius functor (since it is a composition of two Frobenius
functors) and $_A\coring{C}$ is finitely generated projective. On
the other hand, since
$-\otimes_A\coring{C}:\mathcal{M}_A\rightarrow\mathcal{M}%
^{\coring{C}}$ is a left adjoint to $\hom{\coring{C}}{\coring{C}}{-}
:\mathcal{M}^{\coring{C}}\rightarrow\mathcal{M}_A.$ Then
$\hom{\coring{C}}{\coring{C}}{-}$ is a Frobenius functor. Therefore,
$\coring{C}$ is finitely generated projective in
$\mathcal{M}^{\coring{C}},$ and hence in $\mathcal{M}_A.$

\begin{lemma}\label{25}
Let $R$ be the opposite algebra of $^*\coring{C}$.
\begin{enumerate}[(1)]
\item $\coring{C}\in{}^{\coring{C}}\mathcal{M}^{A}$ is
quasi-finite (resp. quasi-finite and injector) as a right
$A$-comodule if and only if $_A\coring{C}$ is finitely generated
projective (resp. $_A\coring{C}$ is finitely generated projective
and $_AR$ is flat). Let $\cohom{A}{\coring{C}}{-}
=-\otimes_AR:\mathcal{M}^{A}\rightarrow\mathcal{M}^{\coring{C}}$ be
the cohom functor.\item If $_A\coring{C}$ is finitely generated
projective and $_AR$ is flat, then
\[
_A\cohom{A}{\coring{C}}{A}_\coring{C}\simeq{}_A%
R_{\coring{C}},
\]
where the right $\coring{C}$-comodule structure of $R$ is defined
as in \cite[Lemma 4.3]{Brzezinski:2002}.
\end{enumerate}
\end{lemma}

\begin{proof}
(1) Straightforward from \cite[Example 4.3]{Gomez:2002}.

(2) From \cite[Lemma 4.3]{Brzezinski:2002}, the forgetful functor
$\mathcal{M}^{\coring{C}}\rightarrow\mathcal{M}%
_A$ is the composition of functors
$\mathcal{M}^{\coring{C}}\rightarrow
\mathcal{M}_{R}\rightarrow\mathcal{M}_A.$ By \cite[Proposition
4.2]{Gomez:2002}, $\cohom{A}{\coring{C}}{-}$ is a left adjoint to
$-\square_{\coring{C}}\coring{C}:\mathcal{M}^{\coring{C}}\rightarrow
\mathcal{M}_A$ which is isomorphic to the forgetful functor
$\mathcal{M}%
^{\coring{C}}\rightarrow\mathcal{M}_A.$ Then
$\cohom{A}{\coring{C}}{-}$ is isomorphic to the composition of
functors
\[
\xymatrix@1{\mathcal{M}_A\ar[r]^{-\otimes
_AR}&\mathcal{M}_{R}\ar[r]& \mathcal{M}^{\coring{C}}.}
\]
In particular, $_A\cohom{A}{\coring{C}}{A}_{\coring{C}}\simeq {}_A(
A\otimes_AR)_{\coring{C}}\simeq{}_AR_{\coring{C}}.$
\end{proof}

\begin{corollary}\label{26}
(\cite[27.10]{Brzezinski/Wisbauer:2003}) \newline Let $\coring{C}$
be an $A$-coring and let $R$ be the opposite algebra of
$^*\coring{C}.$ Then the following statements are equivalent
\begin{enumerate}[(a)]
\item The forgetful functor $F:\mathcal{M}^{\coring{C}%
}\rightarrow\mathcal{M}_A$ is a Frobenius functor;\item
$_A\coring{C}$ is finitely generated projective and
$\coring{C}\simeq R$ as $(A,R)$-bimodules, where $\coring{C}$ is a
right $R$-module via $c.r=c_{(1)}.r(c_{(2)})$, for all
$c\in\coring{C}$ and $r\in R.$
\end{enumerate}
\end{corollary}

\begin{proof}
Straightforward from Theorem \ref{23} and Lemma \ref{25}.
\end{proof}

The following proposition gives sufficient conditions to have that
a morphism of corings is right Frobenius if and only if it is left
Frobenius. Note that it says in particular that the notion of
Frobenius homomorphism of coalgebras over fields (by (b)) or of
rings (by (d)) is independent on the side. Of course, the latter
is well-known.

\begin{proposition}\label{27}
Let $(\varphi,\rho):\coring{C}\rightarrow\coring{D}$ be a
homomorphism of corings such that $_A\coring{C}$, $_B\coring{D}$,
$\coring{C}_A$ and $\coring{D}_B$ are flat. Assume that at least one
of the following holds
\begin{enumerate}[(a)] \item $\coring{C}$ and $\coring{D}$ have a
duality, and $_AB$ and $B_A$ are flat;\item $A$ and $B$ are von
Neumann regular
rings;\item $B\otimes_A\coring{C}$ is coflat in $^{\coring{D}%
}\mathcal{M}$ and $\coring{C}\otimes_AB$ is coflat in $\mathcal{M}%
^{\coring{D}}$ and $_AB $ and $B_A$ are flat;\item $\coring{C}$ and
$\coring{D}$ are coseparable corings.
\end{enumerate} Then the following statements are equivalent
\begin{enumerate}
\item $-\otimes_AB:\mathcal{M}^{\coring{C}%
}\rightarrow\mathcal{M}^{\coring{D}}$ is a Frobenius functor;\item
$B\otimes_A-:{}^{\coring{C}}\mathcal{M}\rightarrow{}^{\coring{D}%
}\mathcal{M}$ is a Frobenius functor.
\end{enumerate}
\end{proposition}

\begin{proof}
Obvious from Theorem \ref{leftrightFrobenius} and Theorem
\ref{leftrightFrobeniusduality}.
\end{proof}

Let us finally show how to derive from our results a remarkable
characterization of the so called \emph{Frobenius corings}.

\begin{corollary}\label{28}
(\cite[27.8]{Brzezinski/Wisbauer:2003})\newline The following
statements are equivalent \begin{enumerate}[(a)] \item the forgetful
functor $\mathcal{M}^{\coring{C}}\rightarrow\mathcal{M}_A$ is
a Frobenius functor;\item the forgetful functor $^{\coring{C}%
}\mathcal{M}\rightarrow{}_A\mathcal{M}$ is a Frobenius functor;\item
there exist an $(A,A) $-bimodule map $\eta:A\rightarrow \coring{C}$
and a $(\coring{C},\coring{C}) $-bicomodule map
$\pi:\coring{C}\otimes_A\coring{C}\rightarrow\coring{C}$ such that
$\pi(\coring{C}\otimes_A\eta) =\coring{C}=\pi(
\eta\otimes_A\coring{C}) .$
\end{enumerate}
\end{corollary}

\begin{proof}
The proof of ``$(1)\Leftrightarrow(4)$'' in Proposition \ref{19a}
for $X=\coring{C}\in{}^{A}\mathcal{M}^{\coring{C}}$ and
$\Lambda=\coring{C}\in{}^{\coring{C}}\mathcal{M}^{A}$ remains
valid for our situation. Finally, notice that the condition (4) in
this case is exactly the condition (c).
\end{proof}

\section{Applications to entwined modules}\label{entwined}

In this section we particularize some our results in Section
\ref{Frobmor} to the category of entwined modules. We adopt the
notations of \cite{Caenepeel/Militaru/Zhu:2002}. We start with
some remarks.

\begin{enumerate}[(1)]

\item  Consider a right-right entwining structure
$(A,C,\psi)\in\mathbb{E}_\bullet^\bullet(k)$
 and a left-left entwining structure $(B,D,\varphi)\in{}_\bullet
 ^\bullet\mathbb{E}(k)$. The category of two-sided entwined modules
$_B^{D}\mathcal{M}(\varphi,\psi)_A^{C}$ defined in \cite[pp.
68--69]{Caenepeel/Militaru/Zhu:2002} is isomorphic to the category
of bicomodules $^{D\otimes B}\mathcal{M}^{A\otimes C}$ over the
associated corings.

\item If $(A,C,\psi)$ and $(A^{\prime },C',\psi')$
belong to $\mathbb{E}_\bullet^{\bullet }(k)$ and are such that
$\psi$ is an isomorphism, then $\psi$ is an isomorphism of corings
(see \cite[Proposition 34]{Caenepeel/Militaru/Zhu:2002}), and
consequently if the coalgebra $C$ is flat as a $k$-module, then the
modules $_A(A\otimes C)$ and $(A\otimes C)_A$ are flat, and \[
^{A\otimes C}\mathcal{M}^{A'\otimes C'}\simeq{}_A%
^{C}\mathcal{M}(\psi^{-1},\psi')_{A'}^{C'}.
\]

\item  Let $(\alpha,\gamma):(A,C,\psi)\rightarrow(A',C',\psi')$ be a 
morphism in
$\mathbb{E}_\bullet^{\bullet }(k)$. We know that
$(\alpha\otimes\gamma,\alpha):A\otimes{C}\rightarrow{A'}\otimes{C'}$
is a morphism of corings. The functor $F$ defined in \cite[Lemma
8]{Caenepeel/Militaru/Zhu:2002} satisfies the commutativity of the
diagram
\[
\xymatrix{\mathcal{M}^{A\otimes
C}\ar[rr]^{-\otimes_AA'}\ar[d]^\simeq & &
\mathcal{M}^{A'\otimes C'}\ar[d]^\simeq \\
\mathcal{M}(\psi)_A^{C}\ar[rr]_{-\otimes_AA'} & &
\mathcal{M}(\psi')_{A'}^{C'},}
\]
where $-\otimes_AA':\mathcal{M}^{A\otimes C}\rightarrow
\mathcal{M}^{A'\otimes C'}$ is the induction functor defined in
\cite[Proposition 5.3]{Gomez:2002}.

\end{enumerate}

We obtain the following result concerning the category of entwined
modules.

\begin{theorem}\label{entwined modules}
Let $(\alpha,\gamma):(A,C,\psi)\rightarrow(A',C',\psi')$ be a
morphism in $\mathbb{E}_\bullet^{\bullet }(k)$, such that $_{k}C$
and $_{k}D$ are flat.
\begin{enumerate} \item The following statements are
equivalent\begin{enumerate}[(a)] \item The functor
$F=-\otimes_AA':\mathcal{M}(\psi)_A^C\rightarrow
\mathcal{M}(\psi')_{A'}^{C'}$ defined in \cite[Lemma
8]{Caenepeel/Militaru/Zhu:2002} is a Frobenius functor;
\item the $ A\otimes{C}-A'\otimes C'$-bicomodule
$(A\otimes C)\otimes_AA'$ is quasi-finite injector as a right
$A'\otimes C'$-comodule and there exists an isomorphism of
$A'\otimes C'-A\otimes C$-bicomodules $\cohom{A'\otimes
C'}{(A\otimes{C})\otimes_AA'}{A'\otimes C'}\simeq
A'\otimes_A(A\otimes{C}).$
\end{enumerate}
Moreover, if $A\otimes C$ and $A'\otimes C'$ are coseparable
corings, then the condition ``injector'' in (b) can be deleted.

 \item If $A$ and $A'$ are QF rings and the module
$(A'\otimes C')_{A'}$ is projective, then the following are
equivalent\begin{enumerate}[(a)] \item The functor
$F=-\otimes_AA':\mathcal{M}(\psi)_A^{C}\rightarrow\mathcal{M}
(\psi')_{A'}^{C'}$ defined in \cite[Lemma
8]{Caenepeel/Militaru/Zhu:2002} is a Frobenius functor; \item the $
A\otimes C-A'\otimes C'$-bicomodule $(A\otimes{C})\otimes_AA'$ is
quasi-finite and injective as a right $A'\otimes C'$-comodule and
there exists an isomorphism of $A'\otimes C'-A\otimes C$-bicomodules
$\cohom{A'\otimes C'}{(A\otimes{C})\otimes_AA'}{A'\otimes C'}\simeq
A'\otimes_A(A\otimes C).$
\end{enumerate}
\end{enumerate}
\end{theorem}

\begin{proof}
Follows from Theorem \ref{23} and Theorem \ref{24}.
\end{proof}

\begin{remark}
Let a right-right entwining structure
$(A,C,\psi)\in\mathbb{E}_\bullet^\bullet(k)$. The coseparability of
the coring ${A}\tensor{}{C}$ is characterized in \cite[Theorem
38(1)]{Caenepeel/Militaru/Zhu:2002} (see also \cite[Corollary
3.6]{Brzezinski:2002}).
\end{remark}

\section{Applications to graded ring theory}\label{graded}

In this section we apply our results in the previous sections to the
category of graded modules by a $G$-set. Let $G$ be a group, $A$ be
a $G$-graded $k$-algebra, and let $X$ be a left $G$-set. The
category $(G,X,A)-gr$ of graded left modules by $X$ is introduced
and studied in \cite{Nastasescu/Raianu/VanOystaeyen:1990}. A study
of the graded ring theory can be found in the recent book
\cite{Nastasescu/VanOystaeyen:2004}. We adopt the notations of
\cite{DelRio:1992} and \cite{Caenepeel/Militaru/Zhu:2002}, and we
begin by giving some useful lemmas.

\subsection{Some useful lemmas}

Let $C=kX$ and $C'=kX'$ be two grouplike coalgebras, where $X$ and
$X'$ are arbitrary sets. We know (see \cite[Example
4]{Caenepeel/Militaru/Zhu:2002}) that the category $\mathcal{M}^{C}$
is isomorphic to the category of $X$-graded modules. Moreover we
have the following:

\begin{lemma}\label{graded1}
For a $k$-module $M$ which is both a $X$-graded and a $X'$-graded
module, the following are equivalent \begin{enumerate}[(a)]
\item $M\in{}^{C'%
}\mathcal{M}^{C}$; \item for every $m\in M,x\in X,x'\in {}X',
\quad {}_{x'}(m_{x}) \in M_{x}$; \item for every $m\in M,x\in
X,x'\in {}X', \quad (_{x'}m)_{x}\in{}_{x'}M$;
\item for every $m\in M,x\in X,x'\in {}X', \quad _{x'}(m_{x}) 
=(_{x'}m)_{x}$.
\end{enumerate}
\end{lemma}

\begin{proof}
Let $M=\bigoplus_{x\in X}M_{x}=\bigoplus_{x'\in X'}{}_{x'}M$. At
first observe that the condition $(a)$ is equivalent to the fact
that the diagram
\[
 \xymatrix{M \ar[rr]^{\rho_M} \ar[d]_{\lambda_M}& & M\otimes kX
\ar[d]^{\lambda_M \otimes kX}
\\kX'\otimes M \ar[rr]_-{kX'\otimes \rho_M} & & kX'\otimes M\otimes
kX\simeq M^{(X'\times X)}}
\]
is commutative.

$(a)\Rightarrow(d)$ Let $m\in M$. We have, $\rho_M(m)=\sum_{x\in
X}m_x \otimes x$ and $\lambda_M(m)=\sum_{x'\in X'}x'\otimes
{_{x'}m}$. From the commutativity of the above diagram,
$$\sum_{x\in X}\sum_{x'\in X'}x'\otimes {_{x'}(m_{x})}\otimes
x=\sum_{x'\in X'}\sum_{x\in X}x'\otimes(_{x'}m) _{x}\otimes x.$$
Hence, $_{x'}(m_{x}) =(_{x'}m) _{x}$ for all $x\in X$, and $x'\in
{}X'$.

$(d)\Rightarrow(b)$ Trivial.

$(b)\Rightarrow(a)$  Let $m\in M$. We have, $$(\lambda_M\otimes
kX)\rho_M(m)=\sum_{x\in X}\sum_{x'\in X'}x'\otimes
{_{x'}(m_{x})}\otimes x,$$ and $$(kX'\otimes
\rho_M)\lambda_M(m)=(kX'\otimes \rho_M)\big(\sum_{x\in X}\sum_{x'\in
X'}x'\otimes {_{x'}(m_{x})}\big)=\sum_{x\in X}\sum_{x'\in
X'}x'\otimes {_{x'}(m_{x})}\otimes x.$$ Hence $(a)$ follows.

$(a)\Leftrightarrow(d)\Leftrightarrow(c)$ Follows by symmetry.
\end{proof}

 Now let $G$ and $G'$ be two groups, $A$ be a $G$-graded
$k$-algebra, $A'$ be a $G'$-graded $k$-algebra, $X$ be a
right $G$-set, and $X'$ be a left $G'$-set. Let $kG$ and
$kG'$ be the canonical Hopf algebras.

 Let $\psi:kX\otimes A\rightarrow A\otimes kX,$ be the map
defined by $[x\otimes a_{g}\mapsto a_{g}\otimes xg]  $, and
$\psi':A'\otimes kX'\rightarrow kX'\otimes A'$, be the map defined
by $[a_{g'}'\otimes x'\mapsto g'x'\otimes a_{g'}'] .$

From \cite[\S 4.6]{Caenepeel/Militaru/Zhu:2002}, we have $(kG,A,kX)
\in\mathbb{DK}_\bullet^\bullet(k)  ,$ $(kG',A',kX')
\in{}_\bullet^\bullet\mathbb{DK}(k) ,$ $(A,kX,\psi)
\in\mathbb{E}_\bullet^\bullet(k) ,$ $(A',kX',\psi')
\in{}_\bullet^\bullet\mathbb{E}( k) ,$ and $\mathcal{M}(kG)_A%
^{kX}\simeq gr-(A,X,G),{}$ $_{A'}^{kX'%
}\mathcal{M}(kG')\simeq(G',X',A')-gr.$

\begin{lemma}\label{graded2}
\begin{enumerate}[(1)] \item Let $M$ be a $k$-module having the 
structure of $X$-graded right
$A$-module and $X'$-graded left $A'$-module.The following are
equivalent
\begin{enumerate}[(a)]\item
$M\in{}_{A'}^{kX'}\mathcal{M}(\psi',\psi)_A^{kX}$; \item the
following conditions hold \begin{enumerate}[(i)]\item $M$ is a
$(A',A) $-bimodule,\item for every $m\in M,x\in X,x'\in X',
\quad_{x'}(m_{x})\in M_{x}$ (or $(_{x'}m)_{x}\in{}_{x'}M$, or
$_{x'}(m_{x}) =(_{x'}m)_{x}$),\item for every $x\in X,$ $M_{x}$ is a
submodule of $_{A'}M$, \item for every $x'\in X',$ $_{x'}M$ is a
submodule of $M_A.$
\end{enumerate}
\end{enumerate}
\item $M\in{}_{A'}^{kX'}\mathcal{M}(\psi',\psi)_A^{kX}$ if and
only if $M$ is an $X'\times X$-graded $A'-A$-bimodule (see \cite[pp.
492--493]{DelRio:1992}).
\end{enumerate}
\end{lemma}

\begin{proof}
 $(1)$ $(a)\Leftrightarrow(b)$ We will use the definition of an object 
in the category
 of two-sided entwined modules 
${}_{A'}^{kX'}\mathcal{M}(\psi',\psi)_A^{kX}$
 (see \cite[pp. 68--69]{Caenepeel/Militaru/Zhu:2002}). By Lemma 
\ref{graded1}, the
condition ``$M\in{}^{kX'}\mathcal{M}^{kX}$'' is equivalent to the
condition $(ii)$. We have moreover that the left $A'$-action on $M$
is $kX$-colinear if and only if for every $x\in X,m\in M_x$,
$\rho_M(a'm)=(a'm)\otimes x$, if and only if for every $x\in X,m\in
M_x$, $a'm\in M_x$, if and only if $(iii)$ holds. By symmetry, the
condition ``the right $A$-action on $M$ is $kX'$-colinear'' is
equivalent to the condition $(iv)$.

$(2)$ The ``if'' part is clear (see \cite[p. 493]{DelRio:1992}). For
the ``only if'' part, put $M=\bigoplus_{(x',x)\in {X'\times
X}}M_{(x',x)}$, where $M_{(x',x)}={}_{x'}(M_{x}) =(_{x'}M)_{x}$. Let
$a'\in A'_{g'}, m\in M$. Since $a'.{}_{x'}(m_x)\in M_x$ and
$a'.{}_{x'}(m_x)\in {}_{g'x'}(M)$, $a'.{}_{x'}(m_x)\in
{}_{a'x'}(M_{x})$. Therefore, $A'_{g'}.{}_{x'}(M_x)\subset
{}_{a'x'}(M_{x})$. By symmetry we obtain
$A'_{g'}M_{(x',x)}A_g\subset M_{g'(x',x)g}=M_{(g'x',xg)}$ ($g\in
G,g'\in G',x\in X, x'\in X'$).
\end{proof}

\subsection{Adjoint pairs and Frobenius pairs of functors between 
categories of graded
modules over $G$-sets}

Throughout this subsection, $G$ and $G'$ will be two groups, $A$
will be a $G$-graded $k$-algebra, $A'$ will be a $G'$-graded
$k$-algebra, $X$ will be a right $G$-set, and $X'$ will be a right
$G'$-set. Let $kG$ and $kG'$ be the canonical Hopf algebras. Let
$\psi:kX\otimes A\rightarrow A\otimes kX$ be the map defined by
$[x\otimes a_{g}\mapsto a_{g}\otimes xg]$. Analogously we define the
map $\psi':kX'\otimes A'\rightarrow A'\otimes kX'$. The coaction and
the counit of the coring $A\otimes kX$ are defined by:
$$\Delta_{A\otimes kX}(a\otimes x)=(a\otimes x)\tensor{A}(1_A\otimes 
x),
\quad \epsilon_{A\otimes kX}(a\otimes x)=a \quad (a\in A,x\in X).$$
\\ An important result is that the corings $A\otimes kX$ and
$A'\otimes kX'$ are coseparable. A proof is clear by using
\cite[Corollary 3.6]{Brzezinski:2002} and \cite[Proposition
101]{Caenepeel/Militaru/Zhu:2002}. A direct proof in the setting
of corings is deduced from \cite[Theorem 3.5]{Brzezinski:2002} by
using the cointegral in the coring $A\otimes kX$ given by
$\delta(a\otimes x\otimes y)=a\delta_{x,y}$ (Kronecker's delta)
for $a\in A,x,y\in X$ (see \cite[26.2]{Brzezinski/Wisbauer:2003}
for the definition of a cointegral in a coseparable coring).

\begin{lemma}\label{graded3}
\begin{enumerate}[(1)]\item $\psi$ is bijective, $(A,kX,\psi^{-1})
\in{}_\bullet^\bullet\mathbb{E}(k) $, and $$_A%
^{kX}\mathcal{M}(kG):={}_A^{kX}\mathcal{M}(\psi^{-1})
\simeq(G,X,A)-gr,$$ where the structure of left $G$-set on $X$ is
given by $g.x=xg^{-1}$ $(g\in G, x\in X)$.
\item Every object of the category $^{A'\otimes
kX'}\mathcal{M}^{A\otimes kX}\simeq{}_{A'}^{kX'}\mathcal{M}((\psi')
^{-1},\psi)_A^{kX}$ can be identified to an $X'\times X$-graded
$A'-A$-bimodule.
\end{enumerate}
\end{lemma}

\begin{proof}
$(1)$ From \cite[Proposition 2]{Caenepeel/Militaru/Zhu:2002}, and
since the Hopf algebra $H=kG$ is cocommutative, $S\circ S=1_{H}$ and
$\overline{S}=S^{-1}=S$ is a twisted antipode of $H.$ Hence $\psi$
is bijective and $(A,kX,\psi^{-1}) \in{}_{\bullet
}^\bullet\mathbb{E}(k)  $ (see \cite[p.
49]{Caenepeel/Militaru/Zhu:2002}), and $_A%
^{kX}\mathcal{M}( kG)  :={}_A^{kX}\mathcal{M}(\psi^{-1})
\simeq(G,X,A)-gr$, since $\psi^{-1}:A\otimes kX\rightarrow kX\otimes
A,$ $[a_{g}\otimes x\mapsto g.x\otimes a_{g}]$, where $g.x=xg^{-1}.$

$(2)$ It follows from $(1)$ and Lemma \ref{graded2}.
\end{proof}

\begin{lemma}\label{graded4}
\begin{enumerate}[(1)]
\item Let $M\in gr-(A,X,G)$, and $N\in (G,X,A)-gr$. We know that $M\in
\mathcal{M}^{A\otimes kX}$ by the coaction: $ \rho_M:M\rightarrow
M\tensor{A}{(A\otimes kX)}$, where $\rho_M(m_x)=m_x\tensor{A}
{(1_A\otimes x)}$, and $N\in {}^{A\otimes kX}\mathcal{M}$ by the
coaction: $\lambda_N:N\rightarrow (A\otimes kX)\tensor{A}{N}$, where
$\lambda_N({}_x n)=(1_A\otimes x) \tensor{A}{{}_x n}$. We have
$M\cotensor{(A\otimes kX)}{N}=M\widehat{\otimes}_AN$, where
$M\widehat{\otimes}_AN$ is the additive subgroup of $M\otimes_A N$
generated by the elements  $m\otimes_A n$ where $x\in X,m\in M_x$,
and $n\in {}_xN$ (see \cite[p. 492]{DelRio:1992}).

\item Let $P$ be an $X\times X'$-graded $A-A'$-bimodule. We have the 
commutative diagram:
\[%
\xymatrix{gr-(A,X,G)\ar[d]^\simeq \ar[rr]^{-\widehat{\otimes}_AP} &
&
gr-( A',X',G') \ar[d]^\simeq  \\
\mathcal{M}^{A\otimes kX}
\ar[rr]^{-\cotensor{(A\otimes kX)}{P}} & & \mathcal{M}^{A'\otimes 
kX'},}
\]
where $-\widehat{\otimes}_AP: gr-(A,X,G)\rightarrow gr-(A',X',G')$
is the functor defined in \cite[p. 493]{DelRio:1992}.
\end{enumerate}
\end{lemma}

\begin{proof}
$(1)$ At first we will show that the right $A$-module $A\otimes kX$
is free. More precisely that the family $\lbrace1_A\otimes x\mid
x\in X\rbrace$ is a right basis of it. It is clear that $a_g\otimes
x=(1_A\otimes xg^{-1})a_g$ for $a_g\in A_g$ and $x\in X$. Now
suppose that $\sum_i (1_A\otimes x_i)a_i=0$, where $x_i\in X, a_i\in
A$. Then  $\sum_i \psi(x_i\otimes a_i)=0$. Since $\psi$ is
bijective, we obtain $\sum_i x_i\otimes a_i=0$. Therefore $a_i=0$
for all $i$. Hence the above mentioned family is a basis of the
right $A$-module $A\otimes kX$. (There is a shorter and indirect
proof of this fact by using that $\psi$ is an isomorphism of
$A$-bimodules, and $\{x\otimes 1_A|x\in X\}$ is a basis of the right
$A$-module $kX\otimes A$).

Finally, suppose that  $\sum_i m_{x_i}\otimes_A {}_{y_i}n\in
M\cotensor{(A\otimes kX)}{N}$. Then
$$ \sum_i m_{x_i}\otimes_A 1_A\otimes x_i\otimes_A {}_{y_i}n=
\sum_i m_{x_i}\otimes_A 1_A\otimes y_i\otimes_A {}_{y_i}n.$$ Hence
(the right $A$-module $A\otimes kX$ is free), $x_i=y_i$ for all $i$,
and $M\cotensor{(A\otimes kX)}{N}\subset M \widehat{\otimes}_AN$.
The other inclusion is obvious.

$(2)$ It suffices to show that the map $M\widehat{\otimes}_AP
\rightarrow M\widehat{\otimes}_AP \tensor{A'}{(A'\otimes kX')}$
defined by $$\Bigg[\sum_{x\in F}m_x \otimes_A {}_x p\longmapsto
\sum_{x'\in X'}\bigg(\sum_{x\in F}m_x \otimes_A{}_x p\bigg)_{x'}
\tensor{A'}{(1_{A'}\otimes x')}\Bigg],$$ where $F$ is a finite
subset of $X$, makes commutative the following diagram
\[
 \xymatrix{0\ar[r] & M\widehat{\otimes}_AP \ar[rr]^{i} \ar[d]
 & & M\tensor{A}{P} \ar[d]
\\ 0\ar[r] & M\widehat{\otimes}_AP\tensor{A'}{(A'\otimes kX')}
\ar[rr]^{i\tensor{A'}{(A'\otimes kX')}} & &
M\tensor{A}{P}\tensor{A'}{(A'\otimes kX')}.}
\]
That is clear since $(\sum_{x\in F}m_x\otimes_A{}_x
p)_{x'}=\sum_{x\in F}m_x\otimes_A{}(_x p)_{x'}=\sum_{x\in F}m_x
\otimes_A{}p_{(x,x')}$, where $p=\sum_{x\in F}{}_x p$.
\end{proof}

Let $\widehat{A}=A\otimes kX$ be the $X\times X$-graded
$A-A$-bimodule associated to the $(A\otimes kX)-(A\otimes kX)$-
bicomodule $A\otimes kX$. It is clear that
$-\widehat{\otimes}_A\widehat{A}= -\cotensor {(A\otimes
kX)}{(A\otimes kX)} \simeq 1_{gr-(A,G,X)}$. The gradings are
$\widehat{A}_x=A\otimes kx$, and ${}_x\widehat{A}=\{\sum_ia_i\otimes
x_i\mid x_ig^{-1}=x, \forall i,\forall g\in G:(a_i)_g\neq 0\}$
($x\in X$). Recall from \cite[Proposition 1.2]{Menini:1993} that for
every $X\times X'$-graded $A-A'$-bimodule $P$, we have an adjunction
$(-\widehat{\otimes}_AP,\h{P_{A'}}{-})$. The unit and the counit of
this adjunction are given respectively by $\eta_M:M\rightarrow
\h{P_{A'}}{M\widehat{\otimes}_AP}$, $\eta_M(m)(p)=\sum_{x\in
X}m_x\tensor{A}{_xp}$ ($m=\sum_{x\in X}m_x\in M,p=\sum_{x\in
X}{}_xp\in P$), and
$\varepsilon_N:\h{P_{A'}}{N}\widehat{\otimes}_AP\rightarrow N$,
$\varepsilon_N(f\tensor{A}{P})=f(p)$ ($f\in \h{P_{A'}}{N}_x,p\in
{}_xP, x\in X$).

\begin{proposition}(\cite[Proposition 1.3, Corollary 1.4]{Menini:1993}) 
\label{graded5}
\begin{enumerate}[(1)]\item The following statements are equivalent for 
a
$k$-linear functor $F:gr-(A,X,G) \rightarrow gr-(A',X',G')$.
\begin{enumerate}[(a)]\item $F$ has a right adjoint; \item $F$ is
right exact and preserves coproducts;
\item $F\simeq -\widehat{\otimes}_AP$ for some $X\times
X'$-graded $A-A'$-bimodule $P$.
\end{enumerate}
\item A $k$-linear functor $G:gr-(A',X',G')\rightarrow gr-(A,X,G)$ has 
a left adjoint
if and only if $G\simeq
\h{P_{A'}}{-}$ for some $X\times X'$-graded $A-A'$-bimodule $P$.
\end{enumerate}
\end{proposition}

\begin{proof}
$(1)$ $(a)\Rightarrow (b)$ Clear. $(b)\Rightarrow (c)$ It follows
from Theorem \ref{3} and Lemma \ref{graded4}. $(c)\Rightarrow (a)$
Obvious from the above mentioned result.

$(2)$ It follows from $(1)$ and \cite[Proposition 1.2]{Menini:1993}.
\end{proof}

\begin{lemma}\label{graded6}
Let $P$ be an $X\times X'$-graded $A-A'$-bimodule.
\begin{enumerate}[(1)]\item $\h{P_{A'}}{-}$ is right exact and 
preserves
direct limits if and only if $_{x}P$ is finitely generated
projective in $\mathcal{M}_{A'}$ for every $x\in X$. \item Suppose
that $_{x}P$ is finitely generated projective in $\mathcal{M}_{A'}$
for every $x\in X$. \begin{enumerate}[(a)] \item For every
$k$-algebra $T$,
$$\xymatrix{\Upsilon_{Z,M}:Z\tensor{T}{\h{P_{A'}}{M}} \ar[r]^-\simeq
& \h{P_{A'}}{Z\tensor{T}{M}}}$$ defined by
$$\Upsilon_{Z,M}(z\tensor{T}{f})(p)=\sum_{x\in
X}z\tensor{T}{f_x(_xp)}$$ $(z\in Z,f=\sum_{x\in X}f_x\in
\h{P_{A'}}{M},p=\sum_{x\in X}{}_{}xp\in P)$, is the natural
isomorphism associated to the functor (see Section
\ref{Frobeniusgeneral})
$$\h{P_{A'}}{-}:gr-(A',X',G')\rightarrow gr-(A,X,G).$$
\item Moreover we have the natural isomorphism
$$\xymatrix{\eta_N:N \widehat{\otimes}_{A'}\h{P_{A'}}{\widehat{A'}}
\ar[r]^-\simeq & \h{P_{A'}}{N}}$$ defined by
$$\eta_N(n_{x'}\tensor{A'}{}_{x'}f)(p)=\sum_{x\in
X}n_{x'}\delta\big((1_{A'}\otimes
x')\tensor{A'}f_{(x',x)}(_xp)\big)$$ $(n_{x'}\in N_{x'},{}_{x'}f\in
{}_{x'}\h{P_{A'}}{\widehat{A'}},p=\sum_{x\in X}{}_xp\in P)$, where
$\delta$ is the cointegral in the coring $A'\otimes kX'$ defined in
the beginning of this subsection. The left grading on
$\h{P_{A'}}{\widehat{A'}}$ is given by
$$_{x'}\h{P_{A'}}{\widehat{A'}}=\bigg\{f\in
\h{P_{A'}}{\widehat{A'}}\mid \Delta_{(A'\otimes
kX')}\big(f(p)\big)=(1_{A'}\otimes x')\tensor{A'}{\sum_{x\in
X}f_x(_{}xp)} \textrm{ for all }p\in P \bigg\}$$ $(x'\in X')$.
\end{enumerate}
\end{enumerate}
\end{lemma}

\begin{proof}
$(1)$ We have
\[
\xymatrix{\h{P_{A'}}{-} \simeq\bigoplus_{x\in
X}\hom{gr-(A',X',G')}{_xP}{-} :gr-(A',X',G') \ar[r] & \mathbf{Ab.}}
\]
Hence, $\h{P_{A'}}{-}$ is right exact and preserves direct limits if
and only if $\hom{gr-(A',X',G')}{_xP}{-}$ is right exact and
preserves direct limits for every $x\in X$ if and only if $_xP$ is
finitely generated projective in $\mathcal{M}_{A'}$ for every $x\in
X$ (by \cite[Proposition V.3.4]{Stenstrom:1975} and the fact that
every finitely generated projective object of a Grothendieck
category is finitely presented).

$(2)$ $(a)$ Let $T$ be a $k$-algebra. Let $M$ be an $X_0\times
X'$-graded $T-A'$-bimodule, where $X_0$ is a singleton, $Z$ be a
right $T$-module, and $x\in X$. We have a sequence of $T$-submodules
$$\h{P_{A'}}{M}_x\leq \h{P_{A'}}{M}\leq \hom{gr-(A',X',G')}{P}{M}\leq
{}_T\hom{A'}{P_{A'}}{M},$$ and the induced structure of left
$T$-module on $\h{P_{A'}}{M}$ is the same structure (see Section
\ref{Frobeniusgeneral}) of left $T$-module associated to the
functor $\h{P_{A'}}{-}$ on it. Moreover, we have the isomorphism
of left $T$--modules $\h{P_{A'}}{M}_x\simeq
\hom{gr-(A',X',G')}{_xP}{M}$. From \cite[Proposition
20.10]{Anderson/Fuller:1992}, for every $x\in X$, there is an
isomorphism
\begin{equation} \label{eta x}
\xymatrix{\eta_x:Z\tensor{T}{\hom{A'}{_xP}{M}} \ar[r]^-\simeq &
\hom{A'}{_xP}{Z\tensor{T}{M}}}
\end{equation}
defined by $\eta_x(z\tensor{T}{\gamma_x}):{}_xp\mapsto
z\tensor{T}\gamma_x(_xp)$ $(z\in Z,\gamma_x\in \hom{A'}{_xP}{M})$.
For every $x\in X$, $\eta_x$ induces an isomorphism
\begin{equation} \label{eta'x}
\xymatrix{\eta'_x:Z\tensor{T}{\hom{gr-(A',X',G')}{_xP}{M}}
\ar[r]^-\simeq & \hom{gr-(A',X',G')}{_xP}{Z\tensor{T}{M}}}.
\end{equation}
Now let us consider the isomorphism
$$\xymatrix{\eta'_{Z,M}:=\bigoplus_{x\in 
X}\eta'_x:Z\tensor{T}{\h{P_{A'}}{M}} \ar[r]^-\simeq
& \h{P_{A'}}{Z\tensor{T}{M}}}.$$ We have
$\eta'_{Z,M}(z\tensor{T}{f})(p)=\sum_{x\in X}z\tensor{T}{f_x(_xp)}$
$(z\in Z,f=\sum_{x\in X}f_x\in \h{P_{A'}}{M},p=\sum_{x\in X}{}_xp\in
P)$. We can verify easily that $\eta'_{Z,M}$ is a morphism of right
$A$-modules. Hence it is a morphism in $gr-(A,X,G)$. It is clear
that $\eta'_{Z,M}$ is natural in $Z$, and $\eta'_{T,M}$ makes
commutative the diagram \eqref{T} (see Section
\ref{Frobeniusgeneral}). Finally, by Mitchell's Theorem
\cite[Theorem 3.6.5]{Mitchell:1965}, $\eta'_{Z,M}=\Upsilon_{Z,M}.$
\\ $(b)$ To prove the first statement it suffices to use the
proof of Theorem \ref{3}, the property $(a)$, and the fact that
every comodule over an $A$-coseparable coring is $A$-relative
injective comodule (for the last see
\cite[26.1]{Brzezinski/Wisbauer:2003}).
\\Finally, we will prove the
last statement. We know that $\lambda_{\h{P_{A'}}{\widehat{A'}}}$ is
defined (\cite{Gomez:2002}) to be the unique $A'$-linear map making
commutative the following diagram
\[
\xymatrix{\h{P_{A'}}{\widehat{A'}}
\ar[rr]^{\lambda_{\h{P_{A'}}{\widehat{A'}}}}
\ar[dr]_{\h{P_{A'}}{\Delta_{(A'\otimes kX')}}} & & \widehat{A'}
\tensor{A'} \h{P_{A'}}{\widehat{A'}}
\ar[dl]^{\Upsilon_{\widehat{A'},\widehat{A'}}}
\\  & \h{P_{A'}}{\widehat{A'}\tensor{A'}\widehat{A'}} &
}.
\]
On the other hand,  for each $x'\in X'$,
$_{x'}\h{P_{A'}}{\widehat{A'}}=\{f\in \h{P_{A'}}{\widehat{A'}}\mid
\lambda_{\h{P_{A'}}{\widehat{A'}}}(f)=(1_{A'}\otimes
x')\tensor{A'}f\}$. Since $\Upsilon_{\widehat{A'},\widehat{A'}}$ is
an isomorphism, $_{x'}\h{P_{A'}}{\widehat{A'}}=\{f\in
\h{P_{A'}}{\widehat{A'}}\mid \Delta_{(A'\otimes
kX')}(f(p))=(1_{A'}\otimes x')\tensor{A'}{\sum_{x\in X}f_x({_xp})}
\textrm{ for all }p\in P \}$.
\end{proof}

\begin{remark}\label{graded7}
Let $\Sigma$ be a finitely generated projective right $A$-module,
and $M$ be a right $A$-module. It easy to verify that the map
$\alpha_{\Sigma,M}:\hom{A}{\Sigma}{M}\rightarrow
M\tensor{A}{\Sigma^*}$ defined by
$\alpha_{\Sigma,M}(\varphi)=\sum_i\varphi(e_i)\tensor{A}{e_i^*}$,
where $\{e_i,e_i^*\}_i$ is a dual basis of $\Sigma$, is an
isomorphism, with $\alpha_{\Sigma,M}^{-1}(m\tensor{A}{f})(u)=mf(u)$
$(m\in M,f\in \Sigma^*,u\in \Sigma)$.
\\ Now we assume that the condition of Lemma \ref{graded6}(2) holds. It 
is easy to verify that for every $x\in X$,
$\eta_x=\alpha^{-1}_{_xP,Z\tensor{T}{M}}\circ(Z\tensor{T}{\alpha_{{}_xP,M}})$,
where $\eta_x$ is the isomorphism \eqref{eta x}. Therefore
$\eta_x^{-1}=(Z\tensor{T}{\alpha^{-1}_{{}_xP,M}})\circ\alpha_{_xP,Z\tensor{T}{M}}.$
Hence $\Upsilon_{Z,M}^{-1}(g)=\sum_{x\in X}\sum_{i\in
I_x}(Z\tensor{T}{\alpha^{-1}_{{}_xP,M}})(g_x(e_{x,i})\tensor{A'
}{e_{x,i}^*})$, where $\{e_{x,i},e_{x,i}^*\}_{i\in I_x}$ is a dual
basis of $_xP$ $(x\in X)$, and $g=\sum_{x\in X}g_x\in
\h{P_{A'}}{Z\tensor{T}{M}}$. Finally, for each $N\in gr-(A',X',G')$,
$$\xymatrix{\eta^{-1}_N:\h{P_{A'}}{N} \ar[r]^-\simeq & N
\widehat{\otimes}_{A'}\h{P_{A'}}{\widehat{A'}}}$$ is
$\eta^{-1}_N=\Upsilon^{-1}_{N,\widehat{A'}}\circ\h{P_{A'}}{\rho_N}$,
and then $$\eta^{-1}_N(\varphi)=\sum_{x\in X}\sum_{i\in
I_x}\sum_{x'\in
X'}(\varphi_x(e_{x,i}))_{x'}\tensor{A'}\psi'(x'\otimes
e_{x,i}^*(-)),$$ where $\varphi=\sum_{x\in X}\varphi_x\in
\h{P_{A'}}{N}$.
\\ In particular, if $X'=G'=\{e'\}$, then 
$$\xymatrix{\eta^{-1}_N:\h{P_{A'}}{N} \ar[r]^-\simeq & N
\otimes_{A'}\h{P_{A'}}{A'}}$$ is defined by
$$\eta^{-1}_N(\varphi)=\sum_{x\in X}\sum_{i\in
I_x}\varphi(e_{x,i})\tensor{A'} e_{x,i}^* \qquad (\varphi\in
\h{P_{A'}}{N}).$$

\end{remark}


\begin{lemma}\label{graded8}
\begin{enumerate}[(1)] \item Let $N\in {}_{A'}\mathcal{M}^{A\otimes 
kX}.$ Then $N$
is quasi-finite as a right $A\otimes kX$-comodule if and only if
$N_x$ is finitely generated projective in $\lmod{A'}$, for every
$x\in X$. In this case, the cohom functor $\cohom{A\otimes
kX}{N}{-}$ is the composite
\[
  \xymatrix{\mathcal{M}^{A\otimes kX}\ar[r]^-\simeq & gr-(A,X,G)
  \ar[r]^-{-\widehat{\otimes}_A P} & \mathcal{M}_{A'}},
\]
where $P$ is the $X\times X'_0$-graded $A-A'$-bimodule
$\h{_{A'}N}{A'}$ with $X'_0$ is a singleton.

\item Now suppose that $N\in{}^{A'\otimes kX'}\mathcal{M}^{A\otimes
kX}$, and $N_x$ is finitely generated projective in $\lmod{A'}$, for
every $x\in X$. Let $\{e_{x,i},e_{x,i}^*\}_{i\in I_x}$ be a dual
basis of $N_x$ $(x\in X)$. Let $\theta:1_{\rcomod{A\otimes
kX}}\rightarrow -\widehat{\otimes}_A P\tensor{A'}N$ be the unit of
the adjunction $(-\widehat{\otimes}_A P,-\tensor{A'}N)$, and let
$M\in gr-(A,X,G)$. The coaction on $\cohom{A\otimes kX}{N}{M}=M
\widehat{\otimes}_A P$:
$$\rho_{M \widehat{\otimes}_A P}:M \widehat{\otimes}_A
P\rightarrow M \widehat{\otimes}_A P\tensor{A'}(A'\otimes kX')$$ is
the unique $A'$-linear map satisfying the condition:
\begin{equation}\label{gradedqf}
\sum_{x\in X}\sum_{i\in I_x}\rho_{M \widehat{\otimes}_A
P}(m_x\tensor{A}e_{x,i}^*)\tensor{A'}e_{x,i}= \sum_{x\in
X}\sum_{i\in I_x}\sum_{x'\in
X'}m_x\tensor{A}e_{x,i}^*\tensor{A'}(1_{A'}\otimes
x')\tensor{A'}(e_{x,i})_{x'},
\end{equation}
for every $m\in M$.
\end{enumerate}
\end{lemma}
\begin{proof}
$(1)$ It follows from Lemma \ref{graded4}(2), that the functor
$F:=-\widehat{\otimes}_{A'}N:\rmod{A'}=gr-(A',X_{0}',G')\to
gr-(A,X,G)$, where $X_{0}'$ is a singleton, is the composite
\[
\xymatrix{\rmod{A'}=gr-(A',X_{0}',G') \ar[r]^-{-\otimes_{A'} N} &
\mathcal{M}^{A\otimes kX} \ar[r]^-\simeq & gr-(A,X,G)}.
\]
Then, $N$ is quasi-finite as a right $A\otimes kX$-comodule if and
only if $F$ has a left adjoint, if and only if (by Corollary
\ref{graded5}(2)) there exists an $X\times X_{0}'$-graded
$A-A'$-bimodule $P$ such that $F\simeq \h{P_{A'}}{-}$, if and only
if (by Lemma \ref{graded6}, \cite[Lemma 3.2(1)]{Gomez:2002} and
Theorem \ref{3}) there exists an $X\times X_{0}'$-graded
$A-A'$-bimodule $P$ such that $_xP$ is finitely generated projective
in $\mathcal{M}_{A'}$ for every $x\in X,$ and $N\simeq
\h{P_{A'}}{A'}$ in $^{A'}\mathcal{M}^{A\otimes kX}.$

Now let us consider the $X\times X'_0$-graded $A-A'$-bimodule
$P:=\h{_{A'}N}{A'}$, and $X'\times X$-graded $A'-A$-bimodule
$M:=\h{P_{A'}}{A'}$, where $X'_0$ is a singleton. We have
$$M=\{f\in P^*|f(_xP)=0 \textrm{ for almost all }x\in X \}$$
$$M_x=\big\{f\in P^*|f(_yP)=0 \textrm{ for all }y\in X-\{x\}\big\}\quad 
(x\in X),$$
where $P^*=\hom{A'}{P_{A'}}{A'_{A'}}$. The structure of
$A'-A$-bimodule on $P^*$ is given by $(fa)(p)=f(ap),
(a'f)(p)=a'f(p)$ $(f\in P^*, a\in A, a'\in A', p\in P)$. We have
$$M\leq (P^*)_A,\quad  M_x\leq M\leq {}_{A'}(P^*), \quad
\textrm{and} \quad M_x\simeq (_xP)^* \textrm{ in }\lmod{A'}\quad
(x\in X).$$

Analogously,
$$P=\{f\in {}^*N|f(N_x)=0 \textrm{ for almost all }x\in X \}$$
$$_xP=\{f\in {}^*N|f(N_y)=0 \textrm{ for all }y\in X-\{x\}\}\quad (x\in 
X),$$
where $^*N=\hom{A'}{_{A'}N}{{}_{A'}A'}$. The structure of
$A-A'$-bimodule on $^*N$ is given by $(af)(n)=f(na),
(fa')(n)=a'f(n)$ $(f\in {}^*N, a\in A, a'\in A', n\in N)$. We have
$$P\leq {}_A(^*N),\quad  _xP\leq P\leq (^*N)_{A'},\quad
\textrm{and} \quad _xP\simeq {}^*(N_x)\textrm{ in } \rmod{A'} \quad
(x\in X).$$ Hence, $N_x$ is finitely generated projective in
$\lmod{A'}$, for every $x\in X$ implies that $_xP$ is finitely
generated projective in $\rmod{A'}$, for every $x\in X$, which
implies that $M_x$ is finitely generated projective in $\lmod{A'}$,
for every $x\in X$.

Finally, let us consider for every $x\in X$ the isomorphism of left
$A'$-modules
$$\xymatrix{H_x:N_x\ar[r]^\simeq & (^*(N_x))^* \ar[r]^\simeq & (_xP)^* 
\ar[r]^\simeq & M_x}.$$
We have $H_x(n_x)(\gamma)=\gamma(n_x)$ if $\gamma\in {}_xP$, and
$H_x(n_x)(\gamma)=0$ if $\gamma\in {}_yP$ and $y\in X-\{x\}$. Set
$$\xymatrix{H=\bigoplus_{x\in X}H_x:N \ar[r]^-\simeq & M}.$$
It can be proved easily that $H$ is a morphism of right $A$-modules.
Hence it is an isomorphism of graded bimodules.

$(2)$ We have $\{e_{x,i}^*,\sigma_x(e_{x,i})\}_{i\in I_x}$ is a dual
basis of $\ldual{(N_x)}\simeq {}_xP$, where $\sigma_x$ is the
evaluation map $(x\in X)$. From the proof of Lemma \ref{graded5} and
Remark \ref{graded7}, the unit of the adjunction
$(-\widehat{\otimes}_A P,-\tensor{A'}N)$ is $\theta_M:M\rightarrow
M\widehat{\otimes}_A P\tensor{A'}N$, $\theta_M(m)=\sum_{x\in
X}\sum_{i\in I_x}m_x\tensor{A}e_{x,i}^*\tensor{A'}e_{x,i}$ $(m\in
M)$. By \cite[4.1]{Gomez:2002}, the coaction on $\cohom{A\otimes
kX}{N}{M}=M\widehat{\otimes}_A P$: $\rho_{M\widehat{\otimes}_A
P}:M\widehat{\otimes}_A P\rightarrow M\widehat{\otimes}_A
P\tensor{A'}(A'\otimes kX')$ is the unique $A'$-linear map
satisfying the commutativity of the diagram
\[
\xymatrix{M\ar[rr]^{\theta_M} \ar[d]^{\theta_M} & &
M\widehat{\otimes}_A P\tensor{A'}N
\ar[d]^{\rho_{M\widehat{\otimes}_A P}\tensor{A'}N}
\\ M\widehat{\otimes}_A
P\tensor{A'}N \ar[rr]^-{M\widehat{\otimes}_A P\tensor{A'}\lambda_N }
& &  M\widehat{\otimes}_A P\tensor{A'}(A'\otimes kX')\tensor{A'}N .}
\]
Finally, the commutativity of the above diagram is equivalent to the
condition \eqref{gradedqf}.
\end{proof}

Now we are in a position to state and prove the main results of this
section which characterize adjoint pairs and Frobenius pairs of
functors between categories of graded modules over $G$-sets.

\begin{theorem}
Let $M$ be an $X\times X'$-graded $A-A'$-bimodule, and $N$ be an 
$X'\times X$-graded
$A'-A$-bimodule. Then the following are equivalent
\begin{enumerate}[(1)]
\item $(-\widehat{\otimes}_AM,-\widehat{\otimes}_{A'}N)$ is an adjoint 
pair;
\item $(N\widehat{\otimes}_A-,M\widehat{\otimes}_{A'}-)$ is an adjoint
pair;
\item $N_x$ is finitely generated projective in $\lmod{A'}$, for every
$x\in X$, and $\widehat{A}\widehat{\otimes}_A \h{_{A'}N}{A'} \simeq
M$ as $X\times X'$-graded $A-A'$-bimodules;
\item
there exist bigraded maps
\[
\psi:\widehat{A}\rightarrow M\widehat{\otimes}_{A'}N \text{ and
 }\omega:N \widehat{\otimes}_A M\rightarrow\widehat{A'},
\]
 such that
\begin{equation}
(\omega\widehat{\otimes}_{A'}N)\circ(N \widehat{\otimes}_A\psi)=
\Lambda\text{ and }(M\widehat{\otimes}_{A'}\omega) \circ(\psi
\widehat{\otimes}_AM)=M.
\end{equation}
\end{enumerate}
In particular, $(-\widehat{\otimes}_AM,-\widehat{\otimes}_{A'}N)$ is
a Frobenius pair if and only if
$(M\widehat{\otimes}_{A'}-,N\widehat{\otimes}_A-)$ is a Frobenius
pair.
\end{theorem}

\begin{proof}
We know that the corings $A\otimes kX$ and $A'\otimes kX'$ are 
coseparable. The Theorem
\ref{19a} achieves the proof.
\end{proof}

\begin{theorem}
Let $G$ and $G'$ be two groups, $A$ be a $G$-graded $k$-algebra, $A'$ 
be a $G'$-graded
$k$-algebra,  $X$ be a right $G$-set,
and $X'$ be a right $G'$-set. For a pair of $k$-linear functors
$F:gr-(A,X,G) \rightarrow gr-(A',X',G')$ and $G:gr-(A',X',G')
\rightarrow gr-(A,X,G)$, the following statements are
equivalent\begin{enumerate}[(a)] \item $(F,G)$ is a Frobenius
pair;\item there exist an $X\times X'$-graded $A-A'$-bimodule $M$,
and an $X'\times X$-graded $A'-A$-bimodule $N$, with the following 
properties
\begin{enumerate}[(1)] \item $_xM$ and $N_x$ is finitely generated 
projective in $\rmod{A'}$ and
$\lmod{A'}$ respectively, for every $x\in X$, and $M_{x'}$ and
$_{x'}N$ is finitely generated projective in $\lmod{A}$ and
$\rmod{A}$ respectively, for every $x'\in X'$,
\item $\widehat{A'}\widehat{\otimes}_{A'}
\h{_AM}{A} \simeq N$ as $X'\times X$-graded $A'-A$-bimodules, and

$\widehat{A}\widehat{\otimes}_A \h{_{A'}N}{A'} \simeq M$ as $X\times
X'$-graded $A-A'$-bimodules,
\item $F\simeq -\widehat{\otimes}_AM$ and
$G\simeq -\widehat{\otimes}_{A'}N$.
\end{enumerate}
\end{enumerate}
\end{theorem}

\begin{proof}
Straightforward from Theorem \ref{leftrightFrobenius}.
\end{proof}

\subsection{When is the induction functor $T^*$ Frobenius?}

Finally, let $f:G\rightarrow G'$ be a morphism of groups, $X$ be a
right $G$-set, $X'$ be a right $G'$-set, $\varphi:X\rightarrow X'$
be a map such that $\varphi(xg) =\varphi(x)f(g) $ for every $g\in
G,$ $x\in X.$ Let $A$ be a $G$-graded $k$-algebra, $A'$ be a
$G'$-graded $k$-algebra, and $\alpha:A\rightarrow A'$ be a morphism
of algebras such that $\alpha( A_{g}) \subset A_{f(g)}'$ for every
$g\in G.$

 We have, $\gamma:kX\rightarrow kX'$ such that $\gamma(x)
=\varphi(x)$ for each $x\in X$, is a morphism of coalgebras, and
$(\alpha ,\gamma) :(A,kX,\psi) \rightarrow(A',kX',\psi') $ is a
morphism in $\mathbb{E}_\bullet^\bullet(k).$

 Let $-\otimes_AA': gr-( A,X,G)\rightarrow gr-(A',X',G')$ be
the functor making commutative the following diagram
\[%
\xymatrix{gr-(A,X,G)\ar[d]^\simeq \ar[rr]^{-\otimes_AA'} & &
gr-( A',X',G') \ar[d]^\simeq  \\
\mathcal{M}(kG)_A^{kX}
\ar[rr]^{-\otimes_AA'%
} & & \mathcal{M}( kG')_{A'%
}^{kX'}%
.}
\]
Let $M\in gr-(A,X,G).$ We have $M\otimes_A A'$ is a right
$A'$-module, and
$$\rho^{r}( m_{x}\otimes _Aa_{g'}')
=(m_{x}\otimes_Aa_{g'%
}')  \otimes\varphi( x)  g'.$$ Therefore, $(M\otimes_AA') _{x'}$ is
the subgroup of $M\otimes_AA'$ spanned by the elements of the form
$m_{x}\otimes_Aa_{g'}'$ where $x\in X,$ $g'\in G',$ $\varphi(x)
g'=x',$ $m_{x}\in M_{x},$ $a_{g'}'\in A_{g'}',$ for every $x'\in
X'.$

 Therefore, the functor $-\otimes_AA':gr-(A,X,G) \rightarrow
gr-(A',X',G')$ is exactly the induction functor $T^*$ defined in
\cite[p. 531]{Menini/Nastasescu:1994}. Hence $T^*$ is a Frobenius
functor if and only if the induction functor
$-\otimes_AA':\mathcal{M}^{A\otimes kX}\rightarrow
\mathcal{M}^{A'\otimes kX'}$ is a Frobenius functor (see Section
\ref{entwined}).

 Moreover, we have the commutativity of the following diagram
\[%
\xymatrix{(G,X,A)-gr \ar[d]^\simeq \ar[rr]^{(T^*)'=A'\otimes_A-}
& & (G',X',A')-gr \ar[d]^\simeq\\
{}_A^{kX}\mathcal{M}(\psi^{-1})
\ar[d]^\simeq\ar[rr]^{A'%
\otimes_A-} & & {}_{A'}^{kX'}\mathcal{M}%
((\psi')^{-1}) \ar[d]^\simeq \\
{}^{kX\otimes A}\mathcal{M} \ar[d]^\simeq
\ar[rr]^{A'\otimes_A-}%
 & & {}^{kX'\otimes A'}\mathcal{M} \ar[d]^\simeq
\\
{}^{A\otimes kX}\mathcal{M} \ar[rr]^{A'\otimes_A-}%
 & & {}^{A'\otimes kX'}\mathcal{M}.}
\]

The following consequence of Theorem \ref{entwined modules} and
\cite[Theorem 2.27]{Menini/Nastasescu:1994} give two different
characterizations when the induction functor $T^*$ is a Frobenius
functor.

\begin{theorem}
The following statements are equivalent\begin{enumerate} [(a)]
\item the functor $T^*:gr-(A,X,G)\rightarrow gr-(A',X',G')$
is a Frobenius functor; \item $(T^*(\widehat{A}))_{x'}$ is finitely
generated projective in $\lmod{A}$, for every $x'\in X'$, and there
exists an isomorphism of $X'\times X$ graded $A'-A$-bimodules
\[
\widehat{A'}\widehat{\otimes}_{A'} \h{_AT^*(\widehat{A})}{A}\simeq
(T^*)'(\widehat{A}).
\]
\end{enumerate}
\end{theorem}

As an immediate consequence of Proposition \ref{27}, we obtain

\begin{proposition}
The following are equivalent\begin{enumerate} [(a)]\item the functor
$T^*:gr-(A,X,G) \rightarrow gr-( A',X',G')$ is a Frobenius functor;
\item the functor $(T^*)':(G,X,A)-gr\rightarrow (G',X',A')-gr$ is
a Frobenius functor.
\end{enumerate}
\end{proposition}

\section*{Acknowledgements}
The results in this paper are part of the author's Ph.D. thesis at the 
University
of Granada, Spain. The author would like to thank his advisor, the Professor
Jos\'e~G\'omez-Torrecillas, for its assistance in
the writing of this paper, and for the fruitful discussions.

The author would also like to thank the Professor Edgar Enochs to have
communicated to him the example of a commutative self-injective ring
which is not coherent. He would like to express his gratitude to the
Professor Tomasz Brzezi\'{n}ski for his assistance, for his
carefully reading of this paper, and for his useful comments on
Section~6.

\end{document}